\documentclass[a4paper,12pt]{article} %izvzela 'twoside' fleqn,

\usepackage{color}

\usepackage{amsfonts, amsmath, amsthm, amssymb}
\usepackage[T1]{fontenc}
\usepackage[cp1250]{inputenc}
\usepackage{xcolor}
\usepackage{graphicx}
\usepackage{amssymb}
\usepackage{amsmath}
\usepackage{mathptmx}
\usepackage{helvet}
\usepackage{courier}
\usepackage{txfonts}
\usepackage{tikz} 
\usetikzlibrary{arrows}
\usepackage{type1cm}

\usepackage{verbatim}

\usepackage{graphicx}
\usepackage{epsfig,amscd,amssymb,amsxtra,amsmath,amsthm}
\usepackage{type1cm}
\usepackage[T1]{fontenc}
\usepackage{graphics}
\usepackage[mathscr]{eucal}
\usepackage[all]{xy}
\usepackage{amsmath,amscd}

\newcommand{\orbit}{\mathcal{O}^{\oplus}}

\newtheorem{theorem}{Theorem}[section]
\newtheorem{proposition}[theorem]{Proposition}
\newtheorem{definition}[theorem]{Definition}
\newtheorem{lemma}[theorem]{Lemma}

\newtheorem{remark}[theorem]{Remark}

\newtheorem{corollary}[theorem]{Corollary}
\newtheorem{problem}[theorem]{Problem}
\newtheorem{observation}[theorem]{Observation}

\newcommand{\diam} {\mathop{\rm diam}\nolimits}

\newcommand{\Cl}  {\mathop{\rm Cl}\nolimits}
\newcommand{\Int}  {\mathop{\rm Int}\nolimits}

\begin{document}

\def\joinrel{\mkern-3mu}
\newcommand{\varproj}{\displaystyle \lim_{\multimapinv\joinrel-\joinrel-}}

\title{A transitive homeomorphism on the Lelek fan}
\author{Iztok Bani\v c, Goran Erceg and Judy Kennedy}
\date{}

\maketitle

\begin{abstract}
Let $X$ be a continuum and let $\varphi:X\rightarrow X$ be a homeomorphism.  {To construct a dynamical system $(X,\varphi)$ with interesting dynamical properties,  the continuum $X$ often needs to have some complicated topological structure. }In this paper, we are interested in one such dynamical property: transitivity.   By  now,  various examples of continua $X$ have been constructed in such a way that the dynamical system $(X,\varphi)$ is transitive.  Mostly, they are examples of continua that are not path-connected,  such as the pseudo-arc or the pseudo-circle,  or they are examples of locally connected continua (and every locally connected continuum is path-connected),  Wa\. zewski's universal dendrite and the Sierpi\' nski carpet are such  examples.   
 
 In this paper, we present an example of a dynamical system $(X,\varphi)$, where $\varphi$ is a homeomorphism on the continuum $X$ and $X$ is a path-connected but not  locally connected continuum.  We construct a transitive homeomorphism on the Lelek fan.  As a by-product,  a non-invertible transitive map on the Lelek fan is also constructed.
\end{abstract}
\-
\\
\noindent
{\it Keywords:} Closed relations; Mahavier products; transitive dynamical systems; transitive homeomorphisms;  fans; Lelek fans\\
\noindent
{\it 2020 Mathematics Subject Classification:} 37B02,37B45,54C60, 54F15,54F17

%%%%%%%%%%%%%%%%%%%%%%%%%%%%%%%%%%%%%%%%%%%%%%%%%%%%%%%%%%%%%%%%%%%%%%%%%%%%%%%%%
%%% I N T R O D U C T I O N S
\section{Introduction}
Let $X$ be a continuum and let $\varphi:X\rightarrow X$ be a homeomorphism on $X$.  To construct a dynamical system $(X,\varphi)$ with interesting dynamical properties (such as minimality, transitivity, non-zero entropy, Li-Yorke or DC2 chaotic systems,  or topological mixing),  the continuum $X$ often needs to have some complicated topological structure, for examples see \cite{barge,jan,jan2,jan3,jan4,HM,chris2,chris3,chris4,chris5,seidler}, where { examples of} such continua may be found.  In this paper, we are interested in transitivity of such  dynamical systems $(X,\varphi)$.  By now,  many such continua $X$ have been constructed.   For example, on a circle, any irrational rotation is a transitive homeomorphism.  It is very challenging  to give  a complete list of such examples. Therefore,  we only present a short list of such  continua.  In the given examples,   references to more examples may be found. 
\begin{enumerate}
\item In \cite{judy}, J. ~Kennedy  and in \cite{minc},  P. ~Minc and W. ~R.~ R. ~Transue constructed independently a transitive homeomorphisms on the pseudo-arc,  which is an example of a continuum that is not path-connected (and, therefore, not locally connected).  P. ~Minc and W. ~R.~ R. ~Transue's construction is based on inverse limits of unit intervals. 
\item  { In \cite{cinc},  J.~\v Cin\v c and P.~Oprocha constructed a non-trivial family of transitive homeomorphisms on the pseudo-arc. This generalizes the  results from \cite{judy} and \cite{minc}. }
\item In \cite{handel},  M.~Handel constructed a transitive homeomorphisms on the pseudo-circle,  which is again an example of a continuum that is not path-connected. 
\item In \cite{chris1},  V. ~Mart\' inez-de-la-Vega,  J.~M.~Mart\' inez-Montejano and  C.~Mouron constructed a transitive homeomorphism on Wa\. zevski's universal dendrite, which is an example of a locally connected continuum.  
\item In \cite{jan4},  J.~ Boro\' nski and P.~Oprocha constructed a transitive homeomorphism on the Sierpi\' nski carpet, which is again an example of a  locally connected continuum. 
\end{enumerate}  
In our paper,  we present an example of a dynamical system $(X,\varphi)$, where $\varphi$ is a homeomorphism on $X$, and the continuum  $X$ is a path-connected but not  locally connected continuum.  We construct a transitive homeomorphism on the Lelek fan. 

In \cite{banic1}, a non-trivial family of fans as  Mahavier products $M_{r,\rho}$ of closed relations $L_{r,\rho}$ on $[0,1]$  was constructed. In particular, the main result  in \cite{banic1}  is a theorem, where the Lelek fan is presented as such a Mahavier product $M_{r,\rho}$.  The closed relation $L_{r,\rho}$ on $[0,1]$ is generated by two line segments in $[0,1]\times [0,1]$. Both line segments contain the origin $(0,0)$ and have positive slopes.  One of the segments extends to the top-boundary and the other to the right side boundary of $[0,1]\times [0,1]$. The Mahavier product $M_{r,\rho}$ is generated by the union $L_{r,\rho}$ of the two line segments. 

In this paper, we first construct a transitive mapping $\sigma_{r,\rho}$ on the Lelek fan $M_{r,\rho}$, which is not a homeomorphism. Then we use $\sigma_{r,\rho}$ to construct a homeomorphism $\sigma$ on the Lelek fan $M=\varprojlim(M_{r,\rho},\sigma_{r,\rho})$.

The Lelek fan was constructed by A.~Lelek in \cite{lelek}. After that, several characterizations of the Lelek fan were presented by J. ~J.~ Charatonik, W. ~J. ~Charatonik and S. ~Miklos in \cite{charatonik2}.  An interesting  property of the Lelek fan $X$ is the fact that the set of its end-points is a dense one-dimensional set in $X$. It is also unique, i.e., it is the only non-degenerate smooth fan with a dense set of end-points.  This was proved independently by W.~D.~Bula and L.~Overseegen  in \cite{oversteegen} and by W. ~Charatonik in  \cite{charatonik}. They proved that any two non-degenerate subcontinua of the Cantor fan with a dense set of endpoints are homeomorphic.

The Lelek fan appears naturally in several dynamical settings and in the theory of continua, see \cite{BK,oversteegen,charatonik, charatonik2,lelek} for more information and more references.  But we have not found any  evidence of examples of transitive mappings being constructed  on the Lelek fan.  However,  recently,  some interesting results were obtained about the Lelek fan or continua with similar properties. For example,   J. ~Boro\' nski, P. ~Minc and S.~ \v Stimac show in \cite{jan}  that for any dendrite $D$,  there is a transitive mapping on $D$.  Note that every dendrite is a locally connected continuum while the Lelek fan is not.  Therefore, the Lelek fan does not fit into this family of continua.  There are also other interesting results presented by L.~C. ~Hoehn and C. ~Mouron in \cite{HM}, where various  interesting transitive mappings that are not homeomorphisms are constructed on the Cantor fan. Moreover,  there are more  interesting results about the Lelek fan, obtained recently  by R.~Hernandez-Gutierrez and L.~C.~Hoehn in \cite{hernandez}, where examples of end-point-rigid smooth fans are presented and studied.  %Also, R. Perez-Marco proved the existence of non-trivial totally invariant connected compacta called hedgehogs near the fixed point of a nonlinearizable germ of holomorphic diffeomorphism. It would be interesting to see, if the Lelek fan may appear as such a hedgehog; see \cite{perez} for more references and  information. 

We proceed as follows. In Section \ref{s1}, the basic definitions and results that are needed later in the paper, are presented. In Section \ref{s2}, we review the construction of the Lelek fan $M_{r,\rho}$ from \cite{banic1} as an infinite Mahavier product of a closed relation $L_{r,\rho}$ on $[0,1]$ and prove some preliminary results that are needed later in Sections \ref{s3} and \ref{s4}.  In Section \ref{s3},  we show that the shift map $\sigma_{r,\rho}$ is a transitive map on the Lelek fan $M_{r,\rho}$.  However, the mapping $\sigma_{r,\rho}$ is not a homeomorphism.  In Section \ref{s4},  we examine in detail the inverse limit $M=\varprojlim(M_{r,\rho},\sigma_{r,\rho})$ and the shift map $\sigma$ on it.  Explicitly, we show that  $M$ is also a Lelek fan.  It follows naturally that  the shift map $\sigma$ is a transitive homeomorphism  on the Lelek fan $M$. 

\section{Definitions and Notation}\label{s1}
The following definitions, notation and well-known results will be needed in the paper.
\begin{definition}
We use \emph{\color{blue} $\mathbb N$} to denote the set of positive integers.
\end{definition}
\begin{definition}
Let $X$ be a metric space and let $A\subseteq X$. We say that $A$ is \emph{\color{blue}dense in $X$} if for any non-empty open set $U$ in $X$, $U\cap A\neq \emptyset$. We also say that $A$ is \emph{\color{blue}nowhere dense in $X$}, if $\Int(\Cl(A))=\emptyset$.
\end{definition}
The following observations are well-known results.
\begin{observation}\label{obi1}
Let $X$ be a  metric space and $A,B\subseteq X$ such that $A$ is dense in $X$
 and $B$ is nowhere dense in $X$. Then $A\setminus B$ is dense in $X$.
 \end{observation}
 \begin{observation}\label{obi2}
Let $X$ be a  metric space and $A,B\subseteq X$ such that $A$ is nowhere dense in $X$
 and $B$ is nowhere dense in $X$. Then $A\cup B$ is nowhere dense in $X$.
 \end{observation}
 \begin{observation}\label{obi3}
Let $X$ be a  metric space, let $A\subseteq X$ and let $U$ be a non-empty open subset of $X$. If $A$ is dense in $X$, then $A\cap U$ is dense in $U$.
 \end{observation}
  \begin{observation}\label{obi4}
Let $X$ be a  compact metric space and for each positive integer $n$, let $U_n$ be a dense open set in $X$.  Then $\bigcap_{n=1}^{\infty}U_n$ is dense in $X$.
 \end{observation}

\begin{definition}
Let $X$ be a metric space, $x\in X$ and $\varepsilon>0$. We use $B(x,\varepsilon)$ to denote the open ball, centered at $x$ with radius $\varepsilon$.
\end{definition}
\begin{definition}\label{nadi}
Let $(X,d)$ be a compact metric space. Then we define \emph{\color{blue}$2^X$} by 
$$
2^{X}=\{A\subseteq X \ | \ A \textup{ is a non-empty closed subset of } X\}.
$$
Let $\varepsilon >0$ and let $A\in 2^X$. Then we define  \emph{\color{blue}$N_d(\varepsilon,A)$} by 
$N_d(\varepsilon,A)=\bigcup_{a\in A}B(a,\varepsilon)$.
Let $A,B\in 2^X$. The function ${\color{blue}H_d}:2^X\times 2^X\rightarrow \mathbb R$, defined by
$$
H_d(A,B)=\inf\{\varepsilon>0 \ | \ A\subseteq N_d(\varepsilon,B), B\subseteq N_d(\varepsilon,A)\},
$$
is called \emph{\color{blue}the Hausdorff metric}. The Hausdorff metric is a metric on the set $2^X$, and the metric space $(2^X,H_d)$ is called \emph{\color{blue}the hyperspace of the space $(X,d)$}. 
\end{definition}
\begin{remark}
Let $(X,d)$ be a compact metric space, let $A$ be a non-empty closed subset of $X$,  and let $(A_n)$ be a sequence of non-empty closed subsets of $X$. When we say $\displaystyle A=\lim_{n\to \infty}A_n$ with respect to the Hausdorff metric, we mean $\displaystyle A=\lim_{n\to \infty}A_n$ in $(2^X,H_d)$. 
\end{remark}
The introduced notation from Definition \ref{nadi} is presented in \cite{nadler}, where more information and references  about hyperspaces may be found.
%\begin{definition}
%Let $X$ and $Y$ be compact metric spaces. 
%A function $F: X\rightarrow 2^Y$ is called \emph{\color{blue}a set-valued function} from $X$ to $Y$. We denote set-valued functions $F: X\rightarrow 2^Y$ by \emph{\color{blue}$F: X\multimap Y$}.
%\end{definition}
%\begin{definition}
%A set-valued function  $F : X\multimap Y$ is  \emph{\color{blue}upper semicontinuous  at a point $x_0\in X$},  if for each
%open set  $U\subseteq Y$ such that $F(x_0)\subseteq U$, there is an open set $V$ in $X$ such that
%\begin{enumerate}
%\item $x_0\in V$ and
%\item for each $x\in V$, $F(x)\subseteq U$. 
%\end{enumerate}  
%The set-valued function  $F : X\multimap Y$ is  \emph{\color{blue}upper semicontinuous},  if it is upper semicontinuous at any point $x\in X$.
%\end{definition}
%\begin{definition}
% \emph{\color{blue}The  graph $\Gamma(F)$ of a set-valued function} $F:X\multimap Y$ is the set of
%all points  $(x,y)\in X\times Y$ such that $y \in F(x)$.  The set-valued function $F$ is surjective, if $\bigcup_{x\in X}F(x)=Y$. 
%\end{definition}
%There is a simple characterization of upper semicontinuous set-valued functions
%(\cite[Proposition 11, p.\ 128]{A} and \cite[Theorem 1.2, p.\ 3]{ingram}):
%
%\begin{theorem}
%\label{th:grafi}  Let $X$ and $Y$ be compact metric spaces and $F:X\multimap Y$ a set-valued function. Then $F$ is upper semicontinuous if and only if its
%graph $\Gamma(F)$ is closed in  $X\times Y$. 
%\end{theorem}
\begin{definition}
Let $X$ be a compact metric space and let $G\subseteq X\times X$ be a relation on $X$. If $G\in 2^{X\times X}$, then we say that $G$ is  \emph{\color{blue}a closed relation on $X$}.  
\end{definition}

%\begin{definition}
%Let $X$  be a set and let $G$ be a relation on $X$.  Then we define  
%$$
%G^{-1}=\{(y,x)\in X\times X \ | \ (x,y)\in G\}
%$$
%to be  \emph{\color{blue}the inverse relation of the relation $G$ on $X$}.
%\end{definition}
\begin{definition}
Let $X$ be a compact metric space and let $G$ be a closed relation on $X$. Then we call
$$
\star_{i=1}^{\infty}G=\Big\{(x_1,x_2,x_3,\ldots )\in \prod_{i=1}^{\infty}X \ | \ \textup{ for each positive integer } i, (x_{i},x_{i+1})\in G\Big\}
$$
 \emph{\color{blue}the  Mahavier product of $G$}.
\end{definition}
%\begin{observation}\label{obsi}
%Let $X$ be a compact metric space, let $f:X\rightarrow X$ be a continuous function. 
%Then 
%$$
%\star_{n=1}^{\infty}\Gamma(f)^{-1}=\varprojlim(X,f).
%$$
%%Also, if $F:X\multimap X$ is an upper semi-continuous function, then  
%%$$
%%\star_{n=1}^{\infty}\Gamma(F)^{-1}=\varproj(X,f).
%%$$
%\end{observation}
\begin{definition}
Let $X$ be a compact metric space and let $G$ be a closed relation on $X$.  The function  
$$
{\color{blue}\sigma} : \star_{n=1}^{\infty}G \rightarrow \star_{n=1}^{\infty}G,
$$
 defined by 
$$
\sigma (x_1,x_2,x_3,x_4,\ldots)=(x_2,x_3,x_4,\ldots)
$$
for each $(x_1,x_2,x_3,x_4,\ldots)\in \star_{n=1}^{\infty}G$, 
is called \emph{\color{blue}the shift map on $\star_{n=1}^{\infty}G$}.  
\end{definition}
\begin{definition}
 \emph{\color{blue}A continuum} is a non-empty compact connected metric space.  \emph{\color{blue}A subcontinuum} is a subspace of a continuum, which is itself a continuum.
 \end{definition}
\begin{definition}
Let $X$ be a continuum. 
\begin{enumerate}
\item The continuum $X$ is \emph{\color{blue}unicoherent}, if for any subcontinua $A$ and $B$ of $X$ such that $X=A\cup B$,  the compactum $A\cap B$ is connected. 
\item The continuum $X$ is \emph{\color{blue}hereditarily unicoherent } provided that each of its subcontinua is unicoherent.
\item The continuum $X$ is a \emph{\color{blue}dendroid}, if it is an arcwise connected, hereditarily unicoherent continuum.
\item Let $X$ be a continuum.  If $X$ is homeomorphic to $[0,1]$, then $X$ is \emph{\color{blue}an arc}.   
\item A point $x$ in an arc $X$ is called \emph{\color{blue}an end-point of the arc  $X$}, if  there is a homeomorphism $\varphi:[0,1]\rightarrow X$ such that $\varphi(0)=x$.
\item Let $X$ be a dendroid.  A point $x\in X$ is called an \emph{\color{blue}end-point of the dendroid $X$}, if for  every arc $A$ in $X$ that contains $x$, $x$ is an end-point of $A$.  The set of all end-points of $X$ will be denoted by $E(X)$. 
\item A continuum $X$ is \emph{\color{blue}a simple triod}, if it is homeomorphic to $([-1,1]\times \{0\})\cup (\{0\}\times [0,1])$.
\item A point $x$ in a simple triod $X$ is called \emph{\color{blue}the top-point} or, briefly, the \emph{\color{blue}top of the simple triod $X$}, if  there is a homeomorphism $\varphi:([-1,1]\times \{0\})\cup (\{0\}\times [0,1])\rightarrow X$ such that $\varphi(0,0)=x$.
\item Let $X$ be a dendroid.  A point $x\in X$ is called \emph{\color{blue}a ramification-point of the dendroid $X$}, if there is a simple triod $T$ in $X$ with the top   $x$.  The set of all ramification-points of $X$ will be denoted by $R(X)$. 
\item The continuum $X$ is \emph{\color{blue}a  fan}, if it is a dendroid with at most one ramification point $v$, which is called the top of the fan $X$ (if it exists).
\item Let $X$ be a fan.   For all points $x$ and $y$ in $X$, we define  \emph{\color{blue}$A[x,y]$} to be the arc in $X$ with end-points $x$ and $y$, if $x\neq y$. If $x=y$, then we define $A[x,y]=\{x\}$.
\item Let $X$ be a fan with the top $v$. We say that that the fan $X$ is \emph{\color{blue}smooth} if for any $x\in X$ and for any sequence $(x_n)$ of points in $X$,
$$
\lim_{n\to \infty}x_n=x \Longrightarrow \lim_{n\to \infty}A[v,x_n]=A[v,x].
$$ 
\item Let $X$ be a fan.  We say that $X$ is \emph{\color{blue}a Cantor fan}, if $X$ is homeomorphic to the continuum $\bigcup_{c\in C}A_c$, where $C\subseteq [0,1]$ is the standard Cantor set and for each $c\in C$, $A_c$ is the  {convex} segment in the plane from $(0,0)$ to $(c,1)$.
\item Let $X$ be a fan.  We say that $X$ is \emph{\color{blue}a Lelek fan}, if it is smooth and $\Cl(E(X))=X$. See Figure \ref{figure2}.
\begin{figure}[h!]
	\centering
		\includegraphics[width=30em]{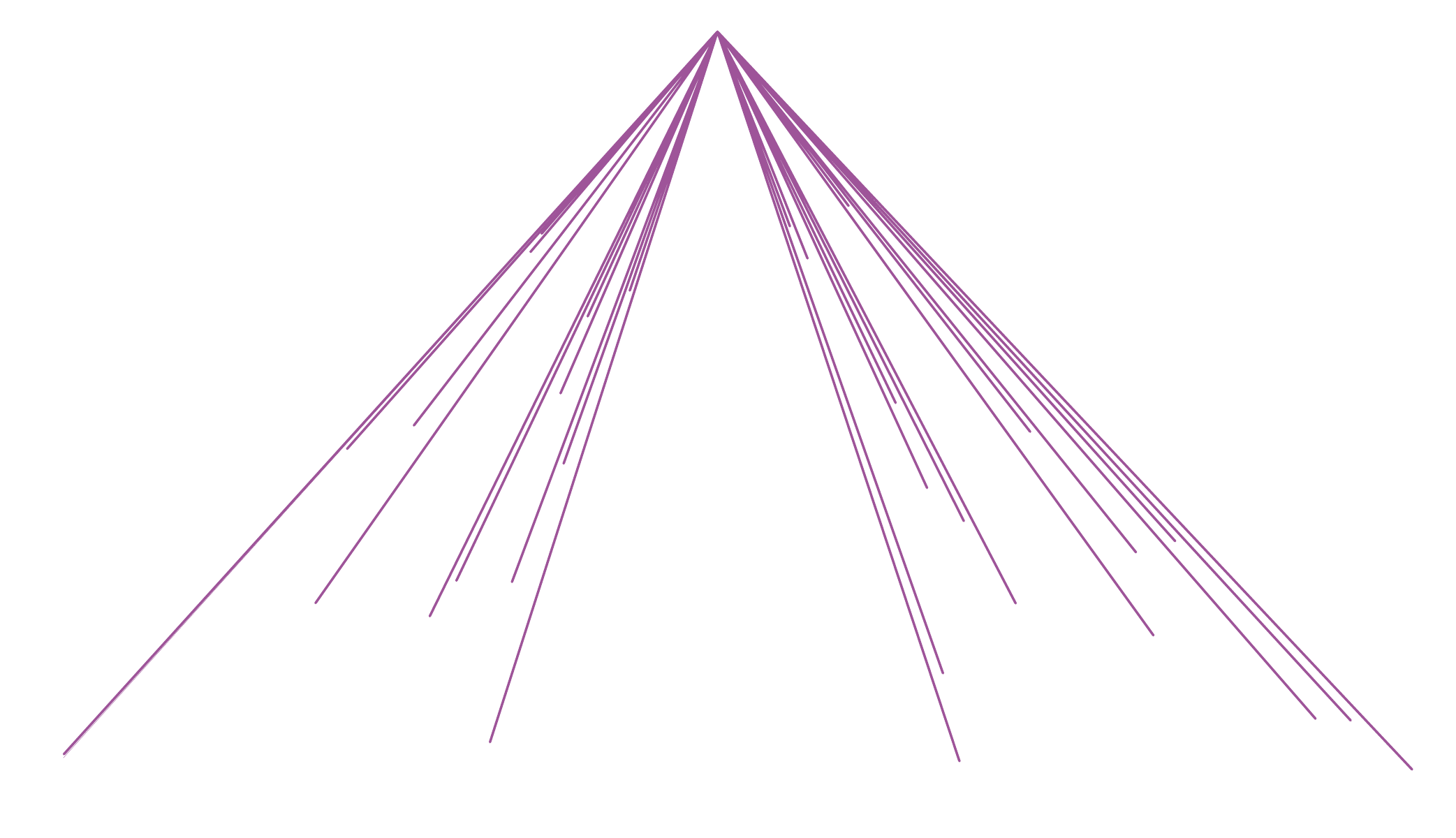}
	\caption{The Lelek fan}
	\label{figure2}
\end{figure}  
\end{enumerate}
See \cite{nadler} for more information about continua, fans and their properties.
\end{definition}

\begin{observation}\label{ffssee}
Let $X$ be the Lelek fan with top $v$.  Note that the set 
$X\setminus(\{v\}\cup E(X))$  is dense in $X$. 
\end{observation}
\begin{definition}
Let $(X,f)$ be a dynamical system.  We say that $(X,f)$ is 
\begin{enumerate}
\item \emph{\color{blue}transitive}, if for all non-empty open sets $U$ and $V$ in $X$,  there is a non-negative integer $n$ such that $f^n(U)\cap V\neq \emptyset$.
\item \emph{\color{blue}dense orbit transitive},   if there is a point $x\in X$ such that its orbit ${\color{blue}\orbit_f(x)}=\{x,f(x),f^2(x),f^3(x),\ldots\}$ is dense in $X$.
\end{enumerate}
\end{definition}
\begin{definition}\label{povezava}
Let $(X,f)$ be a dynamical system. We say that the mapping $f$ is \emph{\color{blue}transitive}, if $(X,f)$ is transitive.
\end{definition}
\begin{observation}\label{isolatedpoints}
It is a well-known fact that if $X$ has no isolated points, then $(X,f)$ is transitive if and only if $(X,f)$ is dense orbit transitive. See \cite{A,KS} for more information about transitive dynamical systems.
\end{observation}
\begin{definition}
Let $(X,f)$ be a dynamical system and let $x_0\in X$ and let $\mathbf x$ be the sequence defined by $\mathbf x=(x_0,f(x_0),f^2(x_0),f^3(x_0),\ldots)$. The set 
$$
{\color{blue}\omega_f(x_0)}=\{x\in X \ | \ \textup{ there is a subsequence of  the sequence } \mathbf x  \textup{ with limit } x\}
$$   
is called \emph{\color{blue}the omega limit set of $x_0$}. 
\end{definition}
\begin{observation}\label{nidani}
It is a well-known fact that $(X,f)$ is transitive if and only if there is a point $x_0\in X$ such that 
$\omega_f(x_0)=X$.  See \cite{A,KS} for more information about omega limit sets of dynamical systems.
\end{observation}
\begin{definition}
Suppose $X$ is a compact metric space.  If $f:X \to X$ is a continuous function, the \emph{\color{blue}inverse limit space} generated by $f$ is the subspace
\begin{equation*}
{\color{blue} \varprojlim(X,f)}=\Big\{(x_{1},x_{2},x_{3},\ldots ) \in \prod_{i=1}^{\infty} X \ | \ 
\text{ for each positive integer } i,x_{i}= f(x_{i+1})\Big\}
\end{equation*}
of the topological product $\prod_{i=1}^{\infty} X$.
\end{definition}
The following theorem is a well-known result, it was proved in  \cite{he} and it can also easily follow  from the results from \cite{li}.  Since the proof is short, we present it here.
\begin{theorem}\label{shifttransitive}
Let $(X,f)$ be a dynamical system. If  $(X,f)$ is transitive and $\sigma :\varprojlim(X,f)\rightarrow \varprojlim(X,f)$ is the shift map on $\varprojlim(X,f)$, then also $(\varprojlim(X,f),\sigma)$ is a transitive dynamical system.  
\end{theorem}
\begin{proof}
Note that since $(X,f)$ is transitive, there is a point $x_0\in X$ such that $\omega_f(x_0)=X$.  Let $x_0$ be such a point.  By \cite[Theorem A, page 96]{li}, for each $\mathbf x=(x_1,x_2,x_3,\ldots)\in \varprojlim (X,f)$ such that $x_1=x_0$, 
$\omega_{\sigma}(\mathbf x)=\varprojlim(\omega_f(x_0),f|_{\omega_f(x_0)})$. 
Since $\omega_f(x_0)=X$, it follows that 
$$
\omega_{\sigma}(\mathbf x)=\varprojlim(\omega_f(x_0),f|_{\omega_f(x_0)})=\varprojlim(X,f).
$$
By Observation \ref{nidani},  $(\varprojlim(X,f),\sigma)$ is a transitive dynamical system.  
\end{proof}
\begin{observation}\label{homeomorphism}
For any dynamical system $(X,f)$, the shift map  $\sigma$ on $\varprojlim(X,f)$ is a homeomorphism from $\varprojlim(X,f)$ to $\varprojlim(X,f)$.
\end{observation}
We will also use the following definition about multiplying open intervals with positive numbers.
\begin{definition}
	Let  $a,b\in \mathbb R$ be such that $a<b$ and let $t>0$. Then we define \emph{\color{blue}the product  of $t$ with the open interval $(a,b)$} as follows: $t\cdot (a,b)=(ta,tb)$. 
\end{definition}

\section{The Lelek fan as an infinite Mahavier product}\label{s2}
Here, we recall the basic construction from \cite{banic1} of the Lelek fan as the infinite Mahavier product $M_{r,\rho}$ of a closed relation on $[0,1]$.  Then, we produce our basic tools that will be used later in Sections \ref{s3} and \ref{s4} when proving our main results.  We begin with the following definitions and notation that are introduced in \cite{banic1}. 

\begin{definition}
	For each $(r,\rho)\in (0,\infty)\times (0,\infty)$, we define the sets \emph{\color{blue}$L_r$}, \emph{\color{blue}$L_{\rho}$} and \emph{\color{blue}$L_{r,\rho}$}  as follows:
	\begin{align*}
	{\color{blue}L_r}&=\{(x,y)\in [0,1]\times [0,1] \ | \ y=rx\},
	\\
	{\color{blue}L_{\rho}}&=\{(x,y)\in [0,1]\times [0,1] \ | \ y=\rho x\},
	\\
	{\color{blue}L_{r,\rho}}&=L_r\cup L_{\rho}.
	\end{align*}
\end{definition}

\begin{definition}
	For each $(r,\rho)\in (0,\infty)\times (0,\infty)$, we define the set \emph{\color{blue}$M_{r,\rho}$}  as follows:
	$$
	{\color{blue}M_{r,\rho}}=\star_{i=1}^{\infty}L_{r,\rho}.
	$$
\end{definition}
\begin{definition}
	Let {$(r,\rho)\in (0,\infty)\times (0,\infty)$}. We say that \emph{\color{blue}$r$ and $\rho$ never connect} or \emph{\color{blue}$(r,\rho)\in \mathcal{NC}$}, if \begin{enumerate}
		\item $r<1$, $\rho>1$ and 
		\item for all integers $k$ and $\ell$,  
		$$
		r^k = \rho^{\ell} \Longleftrightarrow k=\ell=0.
		$$
	\end{enumerate} 
\end{definition}

To prove our main results, we also need the following definitions and results.

\begin{definition}
Let  $(r,\rho)\in \mathcal{NC}$  and let $\mathbf a\in  \{r,\rho\}^{\mathbb N}$. We say that the sequence $\mathbf a$ is \emph{\color{blue}a useful sequence}, if 
$$
{\color{blue}L_{\mathbf a}}{=\{(t,a_1\cdot t,a_2a_1\cdot t,a_3a_2a_1\cdot t,\ldots) \ | \ t\in [0,1]\}\cap M_{r,\rho}}
$$
 is an arc { (i.e.,  $L_{\mathbf a}\neq\{(0,0,0,\ldots)\}$)}.   We also define the set \emph{\color{blue}$\mathcal U_{r,\rho}$} as follows:
$$
\mathcal U_{r,\rho}=\{\mathbf a\in  \{r,\rho\}^{\mathbb N} \ | \ \mathbf a \textup{ is a useful sequence}\}.
$$
\end{definition}
The following lemma is used in Theorem \ref{krajisca}, where a  characterization of end-points of  $M_{r,\rho}$ is obtained.
\begin{lemma}\label{toto}
Let $(r,\rho)\in \mathcal{NC}$,  let $x\in (0,1)$ and let $\mathbf a=(a_1,a_2,a_3,\ldots)\in \{r,\rho\}^{\mathbb N}$ be a sequence such that for each positive integer $n$, $(a_1\cdot a_2\cdot a_3\cdot \ldots \cdot a_n)\cdot x \in [0,1]$. Then
$$
\sup\{(a_1\cdot a_2\cdot a_3\cdot \ldots \cdot a_n)\cdot x \ | \ n \textup{ is a positive integer}\}=1
$$  
if and only if
$$
(x,a_1\cdot x, (a_1\cdot a_2)\cdot x,(a_1\cdot a_2\cdot a_3)\cdot x,\ldots )\in E(M_{r,\rho}).
$$
\end{lemma}
\begin{proof}
The first implication is  \cite[Theorem 4.33]{banic1}. Next, we prove the other implication.  Suppose that 
$$
(x,a_1\cdot x, (a_1\cdot a_2)\cdot x,(a_1\cdot a_2\cdot a_3)\cdot x,\ldots )\in E(M_{r,\rho})
$$
and that $\sup\{(a_1\cdot a_2\cdot a_3\cdot \ldots \cdot a_n)\cdot x \ | \ n \textup{ is a positive integer}\}=s<1$.  Let $t$ be an element of the following open interval:
$$
t\in (s,\min\{1,\sup\{(a_1\cdot a_2\cdot a_3\cdot \ldots \cdot a_n) \ | \ n \textup{ is a positive integer}\}\}).
$$
Also, let 
$$
y=\frac{t}{\sup\{(a_1\cdot a_2\cdot a_3\cdot \ldots \cdot a_n) \ | \ n \textup{ is a positive integer}\}}.
$$
Then 
\begin{align*}
x=&\frac{s}{\sup\{(a_1\cdot a_2\cdot a_3\cdot \ldots \cdot a_n) \ | \ n \textup{ is a positive integer}\}}<
\\
&\frac{t}{\sup\{(a_1\cdot a_2\cdot a_3\cdot \ldots \cdot a_n) \ | \ n \textup{ is a positive integer}\}}=y<
\\
&\frac{\sup\{(a_1\cdot a_2\cdot a_3\cdot \ldots \cdot a_n) \ | \ n \textup{ is a positive integer}\}}{\sup\{(a_1\cdot a_2\cdot a_3\cdot \ldots \cdot a_n) \ | \ n \textup{ is a positive integer}\}}=1.
\end{align*}
and
$$
\sup\{(a_1\cdot a_2\cdot a_3\cdot \ldots \cdot a_n)\cdot y \ | \ n \textup{ is a positive integer}\}=t\leq 1.
$$
Therefore,  the point $(x,a_1\cdot x, (a_1\cdot a_2)\cdot x,(a_1\cdot a_2\cdot a_3)\cdot x,\ldots )$ lies in the interior of the arc $A$ in $M_{r,\rho}$ that contains $(0,0,0,\ldots)$ and $(y,a_1\cdot y, (a_1\cdot a_2)\cdot y,(a_1\cdot a_2\cdot a_3)\cdot y,\ldots )$. Therefore, $(x,a_1\cdot x, (a_1\cdot a_2)\cdot x,(a_1\cdot a_2\cdot a_3)\cdot x,\ldots )$ is not an end-point of $M_{r,\rho}$.
\end{proof}
In Theorem \ref{krajisca}, we give a characterization of end-points of $M_{r,\rho}$. Before we do that, we introduce the following definition that is used in its proof.
\begin{definition}
 For each positive integer $k$,  we use $\pi_k:\prod_{i=1}^{\infty}[0,1]\rightarrow [0,1]$ to denote \emph{\color{blue}the $k$-th standard projection} from $\prod_{i=1}^{\infty}[0,1]$ to $[0,1]$: for each $(x_1,x_2,x_3,\ldots)\in \prod_{i=1}^{\infty}[0,1]$, 
 $$
 \pi_k(x_1,x_2,x_3,\ldots)=x_k.
 $$
\end{definition}
\begin{theorem}\label{krajisca}
Let $(r,\rho)\in \mathcal{NC}$ and let $\mathbf e\in M_{r,\rho}$.  Then $\mathbf e\in E(M_{r,\rho})$ if and only if $\sup\{\pi_n(\mathbf e) \ | \ n \textup{ is a positive integer}\}=1$.
\end{theorem}
\begin{proof}
Let 
$$
\mathbf e=(x,a_1\cdot x, (a_1\cdot a_2)\cdot x,(a_1\cdot a_2\cdot a_3)\cdot x,\ldots )
$$
for some $x\in [0,1]$ and some useful sequence $\mathbf a=(a_1,a_2,a_3,\ldots)$.  Suppose that $\mathbf e\in E(M_{r,\rho})$. Then $x\neq 0$ since $(0,0,0,\ldots)\not \in E(M_{r,\rho})$.  If $x=1$, then obviously, 
$$
\sup\{\pi_n(\mathbf e) \ | \ n \textup{ is a positive integer}\}=1
$$
since $\pi_1(\mathbf 1)=x=1$.  If $x\in (0,1)$, then 
$$
\sup\{\pi_n(\mathbf e) \ | \ n \textup{ is a positive integer}\}=1
$$
by Lemma \ref{toto}. 

To prove the other implication, suppose that $\sup\{\pi_n(\mathbf e) \ | \ n \textup{ is a positive integer}\}=1$.  Obviously, $x\neq 0$ (otherwise, the supremum would equal to $0$).  If $x\in (0,1)$, it follows from Lemma \ref{toto} that $\mathbf e\in E(M_{r,\rho})$.  If $x=1$, obviously, $\mathbf e\in E(M_{r,\rho})$. 
\end{proof}

In \cite{banic1}, the following theorem is the main result.
\begin{theorem}\label{Lelek}
Let $(r,\rho)\in \mathcal{NC}$. Then $M_{r,\rho}$ is a Lelek fan with top $(0,0,0,\ldots)$.
\end{theorem}
\begin{proof}
See the proof of  \cite[Theorem 4.34, page 24]{banic1}.
\end{proof}
\begin{observation} \label{totrebamo}
Note that it follows from Theorems \ref{krajisca} and \ref{Lelek} that for any $\mathbf x=(x_1,x_2,x_3,\ldots)\in M_{r,\rho}$,  
$$
\mathbf x\not \in E(M_{r,\rho}) \Longrightarrow  \text{for each } n\in \mathbb N,  x_n\not\in \{0,1\}.
$$
\end{observation}

The following theorem is also  proved in \cite{banic1}.
%\begin{theorem}\label{dense1}
%	Let  $(r,\rho)\in \mathcal{NC}$.  The set
%	$$
%	\{r^k\cdot \rho^{\ell} \ | \ k \textup{ and } \ell \textup{ are  integers}\}
%	$$
%	is dense in $(0,\infty)$.
%\end{theorem}
\begin{theorem}\label{taoni}
	Let $(r,\rho)\in \mathcal{NC}$.  Then $\{r^k\cdot \rho^{\ell} \ | \ k \textup{ and } \ell \textup{ are  non-negative integers}\}$ is dense in $(0,\infty)$. 
\end{theorem}
\begin{proof}
See the proof of \cite[Theorem 4.23, page 20]{banic1}.
\end{proof}
To prove our main result, we use Proposition \ref{k,l-dense}. Before proving it, we give the following definition. 
\begin{definition}
{ Let $(r,\rho)\in \mathcal{NC}$ and let} $k, \ell \in \mathbb N$. Then we define 
$$
{\color{blue}\mathcal B_{k,\ell}} = \{ r^m \rho^n \ | \  m \geq k, n \geq \ell\}.
$$
\end{definition}
\begin{proposition}\label{k,l-dense}
{ Let $(r,\rho)\in \mathcal{NC}$. Then} for all $k,\ell\in \mathbb N$, $\mathcal B_{k,\ell} $ is dense in $(0, \infty)$.
\end{proposition}
\begin{proof}
Let 
$$
\mathcal B_1(r,\rho) =\{r^k\cdot \rho^{\ell} \ | \ k \textup{ and } \ell \textup{ are  non-negative integers}\}.
$$
First, note that  for all $k, \ell \in \mathbb N$,
$$
\mathcal B_1(r,\rho) \setminus \mathcal B_{k,\ell} = \left( \bigcup_{m=0}^k \{ r^m \rho^n \ | \   n \geq 0 \}	\right) \bigcup \left( \bigcup_{n=0}^\ell \ | \   \{ r^m \rho^n  \ | \   m \geq 0\} \right)
$$
where both, $ \bigcup_{m=0}^k \{ r^m \rho^n \ | \   n \geq 0 \}$ and $\bigcup_{n=0}^\ell \ | \   \{ r^m \rho^n  \ | \   m \geq 0\}$, are nowhere dense in $(0,\infty)$.

Let $k$  and $\ell$ be any positive integers. Then $\mathcal B_1(r,\rho) \setminus \mathcal B_{k,\ell}$ is a union of two nowhere dense sets in $(0,\infty)$ and, therefore, by Observation \ref{obi2}, it is itself a nowhere dense set in $(0,\infty)$. Since the set $\mathcal B_1(r,\rho)$ is dense in $(0,\infty)$, it follows from Observation \ref{obi1}, that $\mathcal B_{k,\ell}$ is dense in $(0, \infty)$.
\end{proof}
\begin{corollary}\label{corkl}
{ Let $(r,\rho)\in \mathcal{NC}$. Then}  for all $k,\ell\in \mathbb N$, $\mathcal B_{k,\ell}  \cap (0,1)$ is dense in $(0,1)$.
\end{corollary}
\begin{proof}
Let $k,\ell\in \mathbb N$. By Proposition \ref{k,l-dense}, $\mathcal B_{k,\ell}$ is dense in $(0,\infty)$. It follows from Observation \ref{obi3} that $\mathcal B_{k,\ell}  \cap (0,1)$ is dense in $(0,1)$.
\end{proof}
Next, we prove Lemma \ref{lemca1}. It will be used in the proof of our main result. To make its statement, we need the following definitions.
\begin{definition}
{ Let $(r,\rho)\in \mathcal{NC}$},  let $\mathbf x=(x_1,x_2,x_3,\ldots)\in M_{r,\rho}$ and let $\mathbf a=(a_1,a_2,a_3,\ldots)\in \{r,\rho\}^{\mathbb N}$ be any sequence of $r$'s and $\rho$'s.  We say that  \emph{\color{blue}the sequence $\mathbf a$ is  associated with $\mathbf x$}, if for each positive integer $n$,
$$
x_{n+1}=a_n\cdot x_n.
$$
\end{definition}
\begin{observation}
{ Let $(r,\rho)\in \mathcal{NC}$} and let 
$\mathbf x\in M_{r,\rho}\setminus\{(0,0,0,\ldots)\}$. Then there is exactly one sequence $\mathbf a\in \{r,\rho\}^{\mathbb N}$ of $r$'s and $\rho$'s such that $\mathbf a$ is associated with $\mathbf x$.
\end{observation}
\begin{definition}
{ Let $(r,\rho)\in \mathcal{NC}$} and let $\mathbf x\in M_{r,\rho}\setminus \{(0,0,0,\ldots)\}$. We use  \emph{\color{blue}
$$
 \mathbf a^{\mathbf x}=(a_1^{\mathbf x},a_2^{\mathbf x},a_3^{\mathbf x},\ldots)
$$
}
 to denote the sequence of $r$'s and $\rho$'s that is associated with $\mathbf x$. 
\end{definition}
To prove the Lemma \ref{lemca1}, we also use a well-known fact, which is stated in the following observation.
\begin{observation}\label{cover}
The following are well-known facts.
\begin{enumerate}
\item Every compact metric space is second countable.
%\item Every subspace of a second countable space is second countable.
\item Let $X$ be a second countable space and $Y\subseteq X$. Then every collection of open subsets of $X$ that covers $Y$ has a countable subcollection that also covers $Y$.
\end{enumerate}
\end{observation}
Finally, we state and prove Lemma \ref{lemca1}. It is used in the proof of Theorem \ref{main}, which is the main result of Section \ref{s3}. 
\begin{lemma}\label{lemca1}
Let $(r,\rho)\in \mathcal{NC}$.  	Then there is a countable collection $\mathcal U$ of open sets in the Hilbert cube $[0,1]^\infty$ such that 
	\begin{enumerate}
		\item\label{1} for each $U \in \mathcal U$, $U\cap M_{r,\rho}\neq \emptyset$;
		\item\label{2} for each $U \in \mathcal U$ there are $n\in \mathbb N$ and open intervals $U_1$, $U_2$, $U_3$,$\ldots$, $U_n$ in $(0,1)$ such that 
		$$
		U=U_1\times \cdots \times U_n \times [0,1]^\infty;
		$$
		\item\label{3} for each $U \in \mathcal U$, for each positive integer $n$ and for all open intervals $U_1$, $U_2$, $U_3$,$\ldots$, $U_n$ in $(0,1)$, it holds that if 
		$$
		U=U_1\times \cdots \times U_n \times [0,1]^\infty, 
		$$
		then there are $a_1, a_2,a_3,\ldots, a_{n-1}\in \{r,\rho\}$ such that 
		\begin{enumerate}
		\item for each $i\in \{1,2,3, \ldots, n-1\}$,
		$$
		U_{i+1}=a_i\cdot U_i, 
		$$
		 \item  for each $\mathbf x \in U \cap M_{r,\rho}$ and for each $i\in \{1,2,3, \ldots, n-1\}$,
		  $$ 
		  a_{i}^{\mathbf x}=a_i;
		  $$
		\end{enumerate}
		\item\label{4} for each $\varepsilon >0$ and for each $\mathbf z=(z_1,z_2,z_3,\ldots) \in M_{r,\rho}$ such that for any positive integer $i$,  $z_i\not\in \{0,1\}$,  there exists $U \in \mathcal U$ such that 
		$$
		\mathbf z \in U \subseteq B(\mathbf z,\varepsilon).
		$$
	\end{enumerate}
\end{lemma}
\begin{proof}
First, let 
$$
M=\{(z_1,z_2,z_3,\ldots)\in M_{r,\rho} \ | \ \text{for each } i\in \mathbb N, z_i\not\in \{0, 1\}\}.
$$
	Next, let  $\mathbf z =(z_1,z_2,z_3,\ldots) \in M$ be any point  and let $\varepsilon  >0$ be arbitrarily chosen.  For each positive integer $i$, let $a_i=a_i^{\mathbf z}$. 
	Let $n$ be a positive integer and let $O_1$,  $O_2$,  $O_3$,  $\ldots$, $O_n$ be open intervals in $(0,1)$ such that 
	$$
	\mathbf z \in  O_1 \times O_2 \times O_3\times \ldots \times O_n \times [0,1]^\infty \subseteq B(\mathbf z,\varepsilon),
	$$
	and let
	$$
	O=O_1 \times O_2 \times O_3\times \ldots \times O_n \times [0,1]^\infty.
	$$
For each $i\in\{1,2,3, \ldots, n-1\}$, let 
$$
\delta_i=\min \{1, \rho z_i \} - rz_i. 
$$
Then $\delta_i>0$ for each $i\in \{1,2,3, \ldots, n-1\}$.  
Also,  let  $\delta >0$ be such that 
$$
\delta < \min \{\delta_i \ | \   i\in \{1,2,3, \ldots, n-1\}\},
$$
and for each $ i \in \mathbb N$, let 
$$
\hat{a}_i=\left\{ \begin{array}{ll}
	r, & \text{ if } a_i=\rho \\
	\rho, & \text{ if } a_i=r
	\end{array} \right..
$$
Note that for each $i\in \{2,3,\ldots,n\}$,
$$
|\hat{a}_iz_i -a_iz_i|=|\hat{a}_i-a_i|z_i=(\rho-r)z_i\geq \delta_i>\delta.
$$
Next, for each $i\in\{1,2,3,\ldots ,n\}$, let $O'_i\subseteq (0,1)$ be an open interval such that
\begin{enumerate}
\item $\diam(O_i')<\delta$ and 
\item $ z_i \in O'_i \subseteq O_i$
\end{enumerate}  
and let 
$$
O'= O'_1 \times O'_2 \times \ldots \times O'_n \times [0,1]^\infty.
$$
Then  
$$
\mathbf z \in   O'_1 \times O'_2 \times \ldots \times O'_n \times [0,1]^\infty \subseteq O_1 \times O_2 \times \ldots \times O_n \times [0,1]^\infty
$$
and the following hold:
\begin{enumerate}
\item $z_1 \in O'_1$, and
\item for each $i\in\{1,2,3,\ldots, n-2\}$, 
$$
a_i\cdot z_i \in O'_{i+1} \text{ but } \hat{a}_iz_i \notin O'_{i+1},
$$ 
since for each $i\in\{1,2,3,\ldots, n-1\}$,  $z_{i+1}=a_{i}\cdot z_{i}$,  $z_{i+1}\in O'_{i+1}$,  $|\hat{a}_{i}\cdot z_{i} -a_{i}\cdot z_{i}|> \delta$ and $\diam(O_{i+1}')<\delta$.
\end{enumerate}
Finally,  let 
	\begin{align*}
	W_n &= O'_n, \\
	W_{n-1} &=  O'_{n-1} \cap  \frac{1}{a_{n-1}} W_{n}=  O'_{n-1} \cap  \frac{1}{a_{n-1}} O'_{n}, \\
	W_{n-2} &=  O'_{n-2} \cap  \frac{1}{a_{n-2}} W_{n-1}=  O'_{n-2} \cap  \frac{1}{a_{n-2}} O'_{n-1}\cap \frac{1}{a_{n-2}a_{n-1}} O'_{n}, \\
	W_{n-3} &=  O'_{n-3} \cap  \frac{1}{a_{n-3}} W_{n-2}=  O'_{n-3} \cap  \frac{1}{a_{n-3}} O'_{n-2}\cap \frac{1}{a_{n-3}a_{n-2}} O'_{n-1}\cap   \frac{1}{a_{n-3}a_{n-2}a_{n-1}} O'_{n},  \\
	&\vdots \\
	W_2 &=O'_{2} \cap  \frac{1}{a_{2}} W_{3}=O'_{2} \cap \frac{1}{a_2} O'_{3} \cap \frac{1}{a_2 a_3} O'_{4} \cap \ldots \cap \frac{1}{a_2a_3a_4 \ldots a_{n-1}} O'_{n}, \\
	W_1 &=O'_{1} \cap  \frac{1}{a_{1}} W_{2}=O'_{1} \cap \frac{1}{a_1} O'_{2} \cap \frac{1}{a_1 a_2} O'_{3} \cap  \frac{1}{a_1 a_2a_3} O'_{4} \cap \ldots \cap \frac{1}{a_1a_2a_3 \ldots a_{n-1}} O'_{n}.
	\end{align*}
Note that for each $i\in \{1,2,3,\ldots,n-1\}$, 
	$$
	a_i W_i \subseteq W_{i+1},
	$$ 	
and for each $i\in \{1,2,3,\ldots, n\}$,  $z_i \in W_i$. Therefore,  
$$
\mathbf z \in W_1 \times W_2 \times \ldots \times W_n \times [0,1]^\infty\subseteq O'\subseteq O\subseteq B(\mathbf z,\varepsilon).
$$
Next,  let $U_1=W_1$ and for each $i\in \{1,2,3,\ldots, n-1\}$, let 
$$
U_{i+1}=a_i\cdot U_i.
$$
Note that for each $i\in\{1,2,3,\ldots, n\}$, $U_i\subseteq W_i$. 
Finally, let 
$$
 U = U_1 \times U_2 \times \ldots \times U_n \times [0,1]^\infty. 
 $$
Recall that at the beginning of the proof, $\mathbf z\in M$ and $\varepsilon>0$ were arbitrarily chosen. So, for any $\mathbf z\in M$ and any $\varepsilon>0$, we construct the open set $U$ as we did in our construction and rename it to $U_{\mathbf z,\varepsilon}$.
Then, let 
$$
\mathcal W =\{   U _{\mathbf z, \varepsilon } \ | \   \mathbf z \in   M  , \varepsilon  >0\}
$$
and let 
$$
D=\bigcap_{i=1}^\infty \left( (0,1)^i \times [0,1]^\infty \right).
$$
Also, for each positive integer $n$, let 
$$
\mathcal W_n=\Big\{W\in \mathcal W \ | \ \diam(W)<\frac{1}{2^n}\Big\}.
$$
Note that for each positive integer $n$,  $\mathcal W_n$ is a collection of open sets in $[0,1]^{\infty}$ such that
$$
D \cap M_{r,\rho}\subseteq \bigcup_{W\in  \mathcal W_n}W.
$$ 
 For each positive integer $n$,  let $ \mathcal V_n$ be a countable subcollection of $\mathcal W_n$ such that 
$$
D \cap M_{r,\rho}\subseteq \bigcup_{V\in \mathcal V_n}V.
$$
Such a subcollection does exist by Observation \ref{cover}.
Finally, let 
$$
\mathcal U = \cup_{n=1}^\infty \mathcal V_n.
$$
Note that $\mathcal U$ is a countable collection  of open sets in the Hilbert cube $[0,1]^\infty$ such that \ref{1},  \ref{2},  \ref{3},  and \ref{4} from Lemma \ref{lemca1} are satisfied. 
\end{proof}

\section{A transitive mapping on the Lelek fan} \label{s3}
Theorem \ref{main} is the main result of this section. It gives a transitive mapping on the Lelek fan. We need the following definitions that are needed in its statement and in its proof.
\begin{definition}
Let $(r,\rho)\in  \mathcal{NC}$. Then we use {\color{blue} $\sigma_{r,\rho}$} to denote the shift map on $M_{r,\rho}$. 
\end{definition}

\begin{definition}
Let $X$ be a set, let $(i_k)$ be a sequence of positive integers, and let $\mathbf x_k=(x_{k,1},x_{k,2},x_{k,3},\ldots ,x_{k,i_k})\in X^{i_k}$ for each positive integer $k$. We define
$$
{\color{blue}\mathbf x_1\oplus \mathbf x_2\oplus \mathbf x_3\oplus\ldots} ={\color{blue}\oplus_{k=1}^{\infty}\mathbf x_k}=(x_{1,1},x_{1,2},x_{1,3},\ldots ,x_{1,i_1},x_{2,1},x_{2,2},x_{2,3},\ldots ,x_{2,i_2},\ldots).
$$
\end{definition}
\begin{theorem}\label{main}
Let $(r,\rho)\in \mathcal{NC}$. Then	$(M_{r,\rho},\sigma_{r,\rho})$ is a transitive dynamical system. 
\end{theorem}
\begin{proof}
Since $M_{r,\rho}$ is a Lelek fan, it is a non-degenerate continuum and therefore, it has no isolated points.  It follows from  Observation \ref{isolatedpoints} that it suffices  to prove that  $(M_{r,\rho},\sigma_{r,\rho})$ is a dense orbit transitive dynamical system. 	So, we prove this theorem by constructing a point $\mathbf z \in M_{r,\rho}$ with a dense orbit in $(M_{r,\rho},\sigma_{r,\rho})$. 
	
	Let $\mathcal U = \{  U_1, U_2, U_3,\ldots \}$
	 be a countable collection of open sets in the Hilbert cube $[0,1]^\infty$ that satisfies the properties \ref{1},  \ref{2},  \ref{3},  and \ref{4} from Lemma \ref{lemca1}. 
For each $i \in \mathbb N$, let
\begin{itemize}
\item $n_i \in \mathbb N$,
\item $a_{i,1}$, $a_{i,2}$, $a_{i,3}$, $\ldots$, $a_{i, n_i-1}\in \{r,\rho\}$, {and}
\item $U_{i,1}$, $U_{i,2}$, $U_{i,3}$, $\ldots$, $U_{i,n_i}$ be open intervals in $(0,1)$,
\end{itemize} 
such that  {(1), (2) and (3) below hold:}
\begin{enumerate}
\item $U_i=U_{i,1} \times U_{i,2}\times U_{i,3} \times \ldots \times U_{i, n_i} \times [0,1]^\infty$,
\item for each $j\in \{1,2,3, \ldots, n_i-1\}$, 
$$
U_{i,j+1}=a_{i,j}\cdot U_{i,j},
$$
\item for each $\mathbf x\in U_i\cap M_{r,\rho}$ and for each $j\in \{1,2,3, \ldots, n_i-1\}$, 
$$
a_{i,j}^{\mathbf x}=a_{i,j}.
$$
\end{enumerate}
	
Next, we construct a point $\mathbf z$ in $M_{r,\rho}$ with a dense orbit in $(M_{r,\rho},\sigma_{r,\rho})$. 

By Corollary \ref{corkl}, $\mathcal{B}_{1,1}\cap (0,1)$ is dense in $(0,1)$. Therefore, 
$$
U_{1,1}\cap(\mathcal{B}_{1,1}\cap (0,1))=U_{1,1}\cap\mathcal{B}_{1,1}\neq \emptyset.
$$
 Let $k_1$ and $\ell_1$ be positive integers such that 
$$
r^{k_1}\cdot \rho^{\ell_1}\in U_{1,1},
$$
and let 
$$
\mathbf z_1=(r,r^2,r^3,\ldots,r^{k_1},r^{k_1}\cdot \rho,r^{k_1}\cdot \rho^{2},r^{k_1}\cdot \rho^{3},\ldots,r^{k_1}\cdot \rho^{\ell_1-1})
$$
and 
$$
\mathbf z_2=(r^{k_1}\cdot \rho^{\ell_1}, a_{1,1}\cdot r^{k_1}\cdot \rho^{\ell_1},a_{1,1}a_{1,2}\cdot r^{k_1}\cdot \rho^{\ell_1},\ldots,a_{1,1}a_{1,2}a_{1,3}\ldots a_{1,n_1-1}\cdot r^{k_1}\cdot \rho^{\ell_1}).
$$
Note that 
$$
\mathbf z_2\in U_{1,1} \times U_{1,2}\times U_{1,3} \times \ldots \times U_{1, n_1}.
$$
Let $k_2'$ and $\ell_2'$ be positive integers such that 
$$
a_{1,1}a_{1,2}a_{1,3}\ldots a_{1,n_1-1}\cdot r^{k_1}\cdot \rho^{\ell_1}=r^{k_2'}\cdot \rho^{\ell_2'}.
$$
By Corollary \ref{corkl}, $\mathcal{B}_{k_2',\ell_2'}\cap (0,1)$ is dense in $(0,1)$. Therefore, 
$$
U_{2,1}\cap(\mathcal{B}_{k_2',\ell_2'}\cap (0,1))=U_{2,1}\cap\mathcal{B}_{k_2',\ell_2'}\neq \emptyset.
$$
 Let $k_2> k_2'$ and $\ell_2> \ell_2'$ be positive integers such that 
$$
r^{k_2}\cdot \rho^{\ell_2}\in U_{2,1},
$$
and let 
$$
\mathbf z_3=(r^{k_2'+1}\cdot \rho^{\ell_2'},r^{k_2'+2}\cdot \rho^{\ell_2'},\ldots,r^{k_2}\cdot \rho^{\ell_2'},r^{k_2}\cdot \rho^{\ell_2'+1},r^{k_2}\cdot \rho^{\ell_2'+2},\ldots,r^{k_2}\cdot \rho^{\ell_2-1})
$$
and 
$$
\mathbf z_4=(r^{k_2}\cdot \rho^{\ell_2}, a_{2,1}\cdot r^{k_2}\cdot \rho^{\ell_2},a_{2,1}a_{2,2}\cdot r^{k_2}\cdot \rho^{\ell_2},\ldots,a_{2,1}a_{2,2}a_{2,3}\ldots a_{2,n_2-1}\cdot r^{k_2}\cdot \rho^{\ell_2}).
$$
Note that 
$$
\mathbf z_4\in U_{2,1} \times U_{2,2}\times U_{2,3} \times \ldots \times U_{2, n_2}.
$$
Let $m$ be any positive integer and suppose that we have already constructed positive integers $k_m'$ and $k_m$,  and $\ell_m'$ and $\ell_m$ such that $k_m>k_m'$, $\ell_m>\ell_m'$ and 
$$
r^{k_m}\cdot \rho^{\ell_m}\in U_{m,1},
$$
and that we have already constructed the points 
$$
\mathbf z_{2m-1}=(r^{k_m'+1}\cdot \rho^{\ell_m'},r^{k_m'+2}\cdot \rho^{\ell_m'},\ldots,r^{k_m}\cdot \rho^{\ell_m'},r^{k_m}\cdot \rho^{\ell_m'+1},r^{k_m}\cdot \rho^{\ell_m'+2},\ldots,r^{k_m}\cdot \rho^{\ell_m-1})
$$
and 
$$
\mathbf z_{2m}=(r^{k_m}\cdot \rho^{\ell_m}, a_{m,1}\cdot r^{k_m}\cdot \rho^{\ell_m},a_{m,1}a_{m,2}\cdot r^{k_m}\cdot \rho^{\ell_m},\ldots,a_{m,1}a_{m,2}a_{m,3}\ldots a_{m,n_m-1}\cdot r^{k_m}\cdot \rho^{\ell_m})
$$
such that
$$
\mathbf z_{2m}\in U_{m,1} \times U_{m,2}\times U_{m,3} \times \ldots \times U_{m, n_m}.
$$
Let $k_{m+1}'$ and $\ell_{m+1}'$ be positive integers such that 
$$
a_{m,1}a_{m,2}a_{m,3}\ldots a_{m,n_m-1}\cdot r^{k_m}\cdot \rho^{\ell_m}=r^{k_{m+1}'}\cdot \rho^{\ell_{m+1}'}.
$$
By Corollary \ref{corkl}, $\mathcal{B}_{k_{m+1}',\ell_{m+1}'}\cap (0,1)$ is dense in $(0,1)$. Therefore, 
$$
U_{{m+1},1}\cap(\mathcal{B}_{k_{m+1}',\ell_{m+1}'}\cap (0,1))=U_{{m+1},1}\cap\mathcal{B}_{k_{m+1}',\ell_{m+1}'}\neq \emptyset.
$$
 Let $k_{m+1}> k_{m+1}'$ and $\ell_{m+1}> \ell_{m+1}'$ be positive integers such that 
$$
r^{k_{m+1}}\cdot \rho^{\ell_{m+1}}\in U_{{m+1},1},
$$
and let 
$$
\mathbf z_{2m+1}=(r^{k_{m+1}'+1}\cdot \rho^{\ell_{m+1}'},r^{k_{m+1}'+2}\cdot \rho^{\ell_{m+1}'},\ldots,r^{k_{m+1}}\cdot \rho^{\ell_{m+1}'},r^{k_{m+1}}\cdot \rho^{\ell_{m+1}'+1},\ldots,r^{k_{m+1}}\cdot \rho^{\ell_{m+1}-1})
$$
and 
$$
\mathbf z_{2m+2}=(r^{k_{m+1}}\cdot \rho^{\ell_{m+1}}, a_{{m+1},1}\cdot r^{k_{m+1}}\cdot \rho^{\ell_{m+1}},\ldots,a_{{m+1},1}a_{{m+1},2}a_{{m+1},3}\ldots a_{{m+1},n_{m+1}-1}\cdot r^{k_{m+1}}\cdot \rho^{\ell_{m+1}}).
$$
Note that 
$$
\mathbf z_{2m+2}\in U_{{m+1},1} \times U_{{m+1},2}\times U_{{m+1},3} \times \ldots \times U_{{m+1}, n_{m+1}}.
$$
Inductively, we have constructed a sequence $(\mathbf z_m)$ such that for each positive integer $m$,
$$
\mathbf z_{2m}\in U_{m,1} \times U_{m,2}\times U_{m,3} \times \ldots \times U_{m, n_m}.$$
Next, let 
$$
\mathbf z=\mathbf z_1\oplus \mathbf z_2\oplus \mathbf z_3\oplus \mathbf z_4\oplus \mathbf z_5\oplus\ldots 
$$
Note that for each positive integer $i$, there is a non-negative integer $\ell$ such that 
$$
\sigma_{r,\rho}^{\ell}(\mathbf z)\in U_i.
$$
To complete the proof, we show that the orbit of $\mathbf z$,
$$
\orbit_{\sigma_{r,\rho}}(\mathbf z)=\{\mathbf z, \sigma_{r,\rho}(\mathbf z), \sigma_{r,\rho}^2(\mathbf z), \sigma_{r,\rho}^3(\mathbf z), \ldots\},
$$
 is dense in $M_{r\rho}$.  To show this,  let $V$ be any open set in the Hilbert cube $[0,1]^{\infty}$ such that $V\cap M_{r,\rho}\neq \emptyset$ and let $U=V\cap M_{r,\rho}$. We prove that
  $$
  U\cap \orbit_{\sigma_{r,\rho}}(\mathbf z)\neq \emptyset.
 $$
 Note that by Observation \ref{ffssee}, the set
 $$
 M_{r,\rho}\setminus (\{(0,0,0,\ldots)\}\cup E(M_{r,\rho}))
 $$
 is dense in $M_{r,\rho}$.
 It follows that there is a point 
 $$
\mathbf x\in  \Big(M_{r,\rho}\setminus (\{(0,0,0,\ldots)\}\cup E(M_{r,\rho}))\Big)\cap U,
 $$
Choose and fix such a point  $\mathbf x=(x_1,x_2,x_3,\ldots)$.  It follows from Observation \ref{totrebamo}, that for each positive integer $i$, $x_i\not\in \{0,1\}$.  Next, let $\varepsilon >0$ be such that $\mathbf x\in B(\mathbf x,\varepsilon)\subseteq V$.   By \ref{4} of Lemma \ref{lemca1},  there exists a positive integer $i$ such that 
		$$
		\mathbf x \in U_i \subseteq B(\mathbf x,\varepsilon).
		$$
Let $i$ be such a positive integer and let $\ell$  be a non-negative integer such that 
$$
\sigma_{r,\rho}^{\ell}(\mathbf z)\in U_i.
$$
Since $U_i \subseteq B(\mathbf x,\varepsilon)\subseteq V$ and $\sigma_{r,\rho}^{\ell}(\mathbf z)\in M_{r,\rho}$, it follows that $\sigma_{r,\rho}^{\ell}(\mathbf z)\in V\cap M_{r,\rho}$ meaning that $\sigma_{r,\rho}^{\ell}(\mathbf z)\in U$.  Therefore,   $  U\cap \orbit_{\sigma_{r,\rho}}(\mathbf z)\neq \emptyset$. This completes the proof.
	\end{proof}
\begin{observation}
Let $(r,\rho)\in \mathcal{NC}$. 
\begin{enumerate}
\item Note that (by Definition \ref{povezava}) the shift map $\sigma_{r,\rho}$  is a transitive map on the Lelek fan $M_{r,\rho}$.  
\item Also, note that  there are points $\mathbf x=(x_1,x_2,x_3,\ldots)\in M_{r,\rho}$ such that $r\cdot x_1\neq \rho\cdot x_1$ and $r\cdot x_1,\rho\cdot x_1\in [0,1]$.  Then 
$$
\sigma_{r,\rho}(r\cdot x_1,x_1,x_2,x_3,\ldots)=\mathbf x \textup{ and } \sigma_{r,\rho}(\rho\cdot x_1,x_1,x_2,x_3,\ldots)=\mathbf x.
$$ 
Therefore, the map $\sigma_{r,\rho}$ is not a homeomorphism.  
	\end{enumerate}
	\end{observation}
\section{A transitive homeomorphism on the Lelek fan}\label{s4}

 In this section, we construct a transitive homeomorphism on the Lelek fan. Let $(r,\rho)\in \mathcal{NC}$.  In this section,  the pair $(r,\rho)$ is fixed.  
\begin{definition}
We use {\color{blue} $M$} to denote the inverse limit
$$
{ {\color{blue}M}}=\varprojlim(M_{r,\rho},\sigma_{r,\rho})
$$
and we use {\color{blue} $\sigma$} to denote the shift map 
$$
{\color{blue}\sigma}:M\rightarrow M
$$
on $M$. 
\end{definition}
\begin{observation}
$M$ is a continuum since it is an inverse limit of continua.
\end{observation}
\begin{observation}
Note that by Observation \ref{homeomorphism},  $\sigma$ is a homeomorphism from $M$ to $M$,  and since $(M_{r,\rho},\sigma_{r,\rho})$ is transitive, it follows from Theorem \ref{shifttransitive},  that $(M,\sigma)$ is a transitive dynamical system. 
\end{observation}
To prove that there is a transitive homeomorphism on the Lelek fan,  we prove that also $M$ is homeomorphic to the  Lelek fan.  
\begin{definition}
We use {\color{blue} $Q$}  to denote the Hilbert cube ${\color{blue}Q}=\prod_{i=1}^{\infty}[0,1]$ and {\color{blue}  $Q^{\infty}$} to denote the countable topological product
$$
{\color{blue}Q^{\infty}}=Q\times Q\times Q\times \ldots=\prod_{i=1}^{\infty}Q.
$$ 
We also use {\color{blue}  $\mathbf O$} to denote the point
$$
{\color{blue}\mathbf O}=((0,0,0,\ldots),(0,0,0,\ldots),(0,0,0,\ldots),\ldots).
$$
\end{definition}
\begin{observation}
Note that $Q^{\infty}$ is homeomorphic to $Q$ and that 
$$
\mathbf O\in M\subseteq Q^{\infty}.
$$
\end{observation}
\begin{definition}
 For each positive integer $k$,  we use ${\color{blue}p_k}:Q^{\infty}\rightarrow Q$ to denote \emph{\color{blue}the $k$-th standard projection} from $Q^{\infty}$ to $Q$: for each $(\mathbf x_1,\mathbf x_2,\mathbf x_3,\ldots)\in Q^{\infty}$, 
 $$
 p_k(\mathbf x_1,\mathbf x_2,\mathbf x_3,\ldots)=\mathbf x_k.
 $$
\end{definition}
\begin{lemma}\label{zejn1}
Let $\mathbf x\in Q^{\infty}$. The following statements are equivalent. 
\begin{enumerate}
\item \label{enica} $\mathbf x\in M\setminus \{\mathbf O\}$. 
\item \label{dvojica} There are
\begin{itemize}
\item a point $x\in (0,1]$,
\item a useful sequence $\mathbf a=(a_1,a_2,a_3,\ldots)\in \mathcal U_{r,\rho}$, and
\item a sequence $\mathbf c=(c_1,c_2,c_3,\ldots)\in \{r,\rho\}^{\mathbb N}$,
\end{itemize} 
such that  { (a), (b) and (c) below hold}
\begin{enumerate}
\item for each positive integer $k$,  $\displaystyle \frac{x}{c_kc_{k-1}c_{k-2}\ldots c_1}<1$, 
\item $p_1(\mathbf x)=(x,a_1\cdot x,a_2a_1\cdot x,a_3a_2a_1\cdot x,\ldots)$, and
\item for each positive integer $k\geq 2$,  
$$
p_k(\mathbf x)=\Big(\frac{x}{c_{k-1}c_{k-2}c_{k-3}\ldots c_3c_2c_1},\ldots ,\frac{x}{c_3c_2c_1},\frac{x}{c_2c_1},\frac{x}{c_1},x,a_1\cdot x,a_2a_1\cdot x,\ldots\Big).
$$
\end{enumerate}
\end{enumerate}
\end{lemma}
\begin{proof}
First, we prove the implication from \ref{enica} to \ref{dvojica}.  For each positive integer $k$, let $\mathbf x^{k}\in M_{r,\rho}$
be such a point that $\mathbf x=(\mathbf x^{1},\mathbf x^{2},\mathbf x^{3},\ldots)$. By Theorem \ref{toto}, for each positive integer $k$, there are a point $x_k\in (0,1]$ and a useful sequence $\mathbf a^k=(a_1^k,a_2^k,a_3^k,\ldots)\in \mathcal U_{r,\rho}$ such that 
$$
\mathbf x^k=(x_k,a_1^k\cdot x_k,a_2^ka_1^k\cdot x_k, a_3^ka_2^ka_1^k\cdot x_k, \ldots).
$$
Let $x=x_1$,  $\mathbf a=\mathbf a^1$ and $\mathbf c=(a_1^2,a_1^3,a_1^4,\ldots)$. We show that 
$$
\mathbf x=\left((x,a_1\cdot x,a_2a_1\cdot x,\ldots),\Big(\frac{x}{c_1},x,a_1\cdot x,a_2a_1\cdot x,\ldots\Big),\Big(\frac{x}{c_2c_1},\frac{x}{c_1},x,a_1\cdot x,a_2a_1\cdot x,\ldots\Big),\ldots\right).
$$
Since $\mathbf x\in M$, it follows that $\mathbf x^1=\sigma_{r,\rho}(\mathbf x^2)$.  Using 
$$
 {\mathbf x^1=}(x,a_1\cdot x,a_2a_1\cdot x, a_3a_2a_1\cdot x, \ldots)=(x_1,a_1^1\cdot x_1,a_2^1a_1^1\cdot x_1, a_3^1a_2^1a_1^1\cdot x_1, \ldots)
$$
and $a_1^2=c_1$ it follows that 
$$
\sigma_{r,\rho}(\mathbf x^2)=(c_1\cdot x_2,a_2^2a_1^2\cdot x_2, a_3^2a_2^2a_1^2\cdot x_2, \ldots)=(x,a_1\cdot x,a_2a_1\cdot x, a_3a_2a_1\cdot x, \ldots)
$$
and, therefore, $x_2=\frac{x}{c_1}$ and
$$
\mathbf x^{2}=\Big(\frac{x}{c_1},x,a_1\cdot x,a_2a_1\cdot x, a_3a_2a_1\cdot x, \ldots\Big)
$$
follows.
Let $k\geq 2$ be a positive integer and suppose that we have already proved that 
$$
\mathbf x^k=(\frac{x}{c_{k-1}c_{k-2}c_{k-3}\ldots c_1},\frac{x}{c_{k-2}c_{k-3}\ldots c_1},\ldots,\frac{x}{c_1},x,a_1\cdot x,a_2a_1\cdot x, a_3a_2a_1\cdot x, \ldots).
$$
Since $\mathbf x\in M$, it follows that $\mathbf x^k=\sigma_{r,\rho}(\mathbf x^{k+1})$.  Using 
$a_1^{k+1}=c_k$ it follows that 
$$
\sigma_{r,\rho}(\mathbf x^{k+1})=(c_k\cdot x_{k+1},a_2^{k+1}a_1^{k+1}\cdot x_{k+1}, a_3^{k+1}a_2^{k+1}a_1^{k+1}\cdot x_{k+1}, \ldots)=\mathbf x^k
$$
and, therefore, $\displaystyle x_{k+1}=\frac{x}{c_kc_{k-1}c_{k-2}c_{k-3}\ldots c_1}$ and
$$
\mathbf x^{k+1}=(\frac{x}{c_kc_{k-1}c_{k-2}\ldots c_1},\frac{x}{c_{k-1}c_{k-2}\ldots c_1},\ldots,\frac{x}{c_1},x,a_1\cdot x,a_2a_1\cdot x, a_3a_2a_1\cdot x, \ldots)
$$
follows. Obviously, for each positive integer $k$,  $\displaystyle \frac{x}{c_kc_{k-1}c_{k-2}\ldots c_1}<1$ (since $\mathbf x\in M$). 

The implication from \ref{dvojica} to \ref{enica} is obvious. 
\end{proof}
\begin{definition}
For all positive integers $m$ and $n$, we use {\color{blue} $p_{m}^{n}$} to denote the projection 
${\color{blue}p_{m}^{n}}:Q^{\infty}\rightarrow [0,1]$
defined by
$$
p_{m}^{n}\Big((x_1^1,x_2^1,x_3^1,\ldots),(x_1^2,x_2^2,x_3^2,\ldots),(x_1^3,x_2^3,x_3^3,\ldots),\ldots\Big)=x_m^n
$$
for all  $\Big((x_1^1,x_2^1,x_3^1,\ldots),(x_1^2,x_2^2,x_3^2,\ldots),(x_1^3,x_2^3,x_3^3,\ldots),\ldots\Big)\in Q^{\infty}$.  
\end{definition}
\begin{observation}
For each $\mathbf x\in M\setminus \{\mathbf O\}$,  there are 
\begin{itemize}
\item the unique  point $x\in (0,1]$; namely, $x=p_1^1(\mathbf x)$,
\item the unique  useful sequence $\mathbf a=(a_1,a_2,a_3,\ldots)\in \mathcal U_{r,\rho}$, and
\item the unique  sequence $\mathbf c=(c_1,c_2,c_3,\ldots)\in \{r,\rho\}^{\mathbb N}$,
\end{itemize} 
such that  
$$
\mathbf x=\left((x,a_1\cdot x,a_2a_1\cdot x,\ldots),\Big(\frac{x}{c_1},x,a_1\cdot x,a_2a_1\cdot x,\ldots\Big),\Big(\frac{x}{c_2c_1},\frac{x}{c_1},x,a_1\cdot x,a_2a_1\cdot x,\ldots\Big),\ldots\right).
$$
\end{observation}

\begin{definition}
For each $\mathbf x\in M\setminus \{\mathbf O\}$,  we use {\color{blue} $\mathbf a(\mathbf x)$} and  {\color{blue}  $\mathbf c(\mathbf x)$} to denote the useful sequence 
$$
{\color{blue}\mathbf a(\mathbf x)}=(a_1,a_2,a_3,\ldots)\in \mathcal U_{r,\rho}
$$
and the sequence 
$$
{\color{blue}\mathbf c(\mathbf x)}=(c_1,c_2,c_3,\ldots)\in \{r,\rho\}^{\mathbb N}
$$
such that 
$$
\mathbf x=\left((x,a_1\cdot x,a_2a_1\cdot x,\ldots),\Big(\frac{x}{c_1},x,a_1\cdot x,a_2a_1\cdot x,\ldots\Big),\Big(\frac{x}{c_2c_1},\frac{x}{c_1},x,a_1\cdot x,a_2a_1\cdot x,\ldots\Big),\ldots\right).
$$
%where $x=p_1^1(\mathbf x)$.
\end{definition}
\begin{definition}
Let $\mathbf a\in \mathcal U_{r,\rho}$ and let $\mathbf c\in \{r,\rho\}^{\mathbb N}$.  We say that the pair ${\color{blue}(\mathbf a,\mathbf c)}$ is \emph{\color{blue}a useful pair}, if there is $\mathbf x\in M\setminus \{\mathbf O\}$ such that $\mathbf a=\mathbf a(\mathbf x)$ and $\mathbf c=\mathbf c(\mathbf x)$. We use {\color{blue}$\mathcal P_{r,\rho}$} to denote 
$$
{\color{blue}\mathcal P_{r,\rho}}=\{(\mathbf a,\mathbf c) \ | \ (\mathbf a,\mathbf c) \textup{ is a useful pair}\}.
$$
\end{definition}
\begin{definition}
For each useful pair $ (\mathbf a,\mathbf c)\in \mathcal P_{r,\rho}$,  {where} $\mathbf a=(a_1,a_2,a_3,\ldots)$ and $\mathbf c=(c_1,c_2,c_3,\ldots)$, for each $x\in [0,1]$, and for each positive integer $k$, we define {\color{blue} $K^{(\mathbf a,\mathbf c)}_k(x)$} as follows:
\begin{enumerate}
\item ${\color{blue}K^{(\mathbf a,\mathbf c)}_1(x)}=(x,a_1\cdot x,a_2a_1\cdot x,a_3a_2a_1\cdot x,\ldots)$ and 
\item for each positive integer $k\geq 2$, 
$$
{\color{blue}K^{(\mathbf a,\mathbf c)}_k(x)}=\Big(\frac{x}{c_{k-1}c_{k-2}c_{k-3}\ldots c_3c_2c_1},\ldots ,\frac{x}{c_3c_2c_1},\frac{x}{c_2c_1},\frac{x}{c_1},x,a_1\cdot x,a_2a_1\cdot x,\ldots\Big).
$$
\end{enumerate}
\end{definition}
\begin{definition}
For each useful pair $ (\mathbf a,\mathbf c)\in \mathcal P_{r,\rho}$,  {where} $\mathbf a=(a_1,a_2,a_3,\ldots)$ and $\mathbf c=(c_1,c_2,c_3,\ldots)$, we define {\color{blue}  $L_{(\mathbf a,\mathbf c)}$} to be the set 
$$
{\color{blue}L_{(\mathbf a,\mathbf c)}}=\left\{\left(K^{(\mathbf a,\mathbf c)}_1(x),K^{(\mathbf a,\mathbf c)}_2(x),K^{(\mathbf a,\mathbf c)}_3(x),\ldots\right) \ | \  x\in [0,1]\right\}\cap M.
$$
\end{definition}
\begin{observation}
Note that 
\begin{enumerate}
\item for each useful pair $ (\mathbf a,\mathbf c)\in \mathcal P_{r,\rho}$,   $L_{(\mathbf a,\mathbf c)}$ is a {convex} segment in $M$ with one end-point being $\mathbf O$,
\item $\displaystyle M=\bigcup_{ (\mathbf a,\mathbf c)\in \mathcal P_{r,\rho}}L_{(\mathbf a,\mathbf c)}$.
\end{enumerate}
\end{observation}

\begin{theorem}\label{fencic}
$M$ is a smooth fan.
\end{theorem}
\begin{proof}
For all $\mathbf a=(a_1,a_2,a_3,\ldots),\mathbf c=(c_1,c_2,c_3,\ldots)\in \{r,\rho\}^{\mathbb N}$, 
let 
$$
{\color{blue}A_{\mathbf a,\mathbf c}}=\Big\{{\Big(K_1^{\mathbf a,\mathbf c}(t),K_2^{\mathbf a,\mathbf c}(t),K_3^{\mathbf a,\mathbf c}(t),\ldots\Big)} \ | \ t\in [0,1]\Big\}
$$
and let 
$$
{\color{blue}F}=\bigcup_{(\mathbf a,\mathbf c)\in \{r,\rho\}^{\mathbb N}\times \{r,\rho\}^{\mathbb N}}A_{\mathbf a,\mathbf c}.
$$
Note that $M\subseteq F$ and that each $A_{\mathbf a,\mathbf c}$ is a convex arc in $P${, where}
$$
{\color{blue}P}=\prod_{j=1}^{\infty}\left(\Big([0,\frac{1}{r^{j-1}}]\times [0,\frac{1}{r^{j-2}}]\times[0,\frac{1}{r^{j-3}}]\times\ldots \times [0,\frac{1}{r^2}]\times [0,\frac{1}{r}]\times [0,1]\Big)\times \prod_{i=1}^{\infty}[0,\rho^{i}]\right).
$$
Therefore, $M\subseteq F\subseteq P$. 

For each positive integer $j$, let 
$$
{\color{blue}Q_j}=\Big([0,\frac{1}{r^{j-1}}]\times [0,\frac{1}{r^{j-2}}]\times[0,\frac{1}{r^{j-3}}]\times\ldots \times [0,\frac{1}{r^2}]\times [0,\frac{1}{r}]\times [0,1]\Big)\times \prod_{i=1}^{\infty}[0,\rho^{i}].
$$
Note that each $Q_j$ is homeomorphic to the Hilbert cube, and since $P=\prod_{j=1}^{\infty}Q_j$, also $P$ is homeomorphic to the Hilbert cube.  For each positive integer $j$,  we define closed intervals ${\color{blue}I_1^j}$, ${\color{blue}I_2^j}$, ${\color{blue}I_3^j}$, $\ldots$ with the following identity:
$$
\prod_{k=1}^{\infty}I_k^j=\Big([0,\frac{1}{r^{j-1}}]\times [0,\frac{1}{r^{j-2}}]\times[0,\frac{1}{r^{j-3}}]\times\ldots \times [0,\frac{1}{r^2}]\times [0,\frac{1}{r}]\times [0,1]\Big)\times \prod_{i=1}^{\infty}[0,\rho^{i}],
$$ 
explicitly,
$$
{\color{blue}I_1^j}=[0,\frac{1}{r^{j-1}}]\textup{,  }{\color{blue}I_2^j}=[0,\frac{1}{r^{j-2}}], \ldots, {\color{blue}I_{j-1}^j}=[0,\frac{1}{r}]\textup{,  }{\color{blue}I_j^j}=[0,1]\textup{,  }{\color{blue}I_{j+1}^j}=[0,\rho]\textup{,  }{\color{blue}I_{j+2}^j}=[0,\rho^2],\ldots
$$
We use ${\color{blue}\sigma_P}:P\rightarrow P$ to denote the shift map on $P$, defined by
$$
\sigma_P\Big((t_1^1,t_2^1,t_3^1,\ldots),(t_1^2,t_2^2,t_3^2,\ldots),(t_1^3,t_2^3,t_3^3,\ldots),\ldots\Big)=\Big((t_1^2,t_2^2,t_3^2,\ldots),(t_1^3,t_2^3,t_3^3,\ldots),\ldots\Big)
$$
for all $\Big((t_1^1,t_2^1,t_3^1,\ldots),(t_1^2,t_2^2,t_3^2,\ldots),(t_1^3,t_2^3,t_3^3,\ldots),\ldots\Big)\in P$, and for each positive integer $j$, we use ${\color{blue}\sigma_j}$ to denote the function $\sigma_j:Q_{j+1}\rightarrow Q_j$, defined by 
$$ 
\sigma_j(x_1,x_2,x_3,x_4,\ldots)=(x_2,x_3,x_4,\ldots)
$$
for each $(x_1,x_2,x_3,x_4,\ldots)\in Q_{j+1}$.

For each positive integer $j$, we use 
$$
{\color{blue}p_j}:P\rightarrow Q_j
$$
to denote the standard projection defined by 
$$
p_j\Big((t_1^1,t_2^1,t_3^1,\ldots),(t_1^2,t_2^2,t_3^2,\ldots),(t_1^3,t_2^3,t_3^3,\ldots),\ldots\Big)=(t_1^j,t_2^j,t_3^j,\ldots)
$$
for each $\Big((t_1^1,t_2^1,t_3^1,\ldots),(t_1^2,t_2^2,t_3^2,\ldots),(t_1^3,t_2^3,t_3^3,\ldots),\ldots\Big)\in P$.
For all positive integers $j$ and $k$, we use 
$$
{\color{blue}\pi^j_k}: Q_j\rightarrow I_k^j
$$
to denote the standard projection defined by 
$$
\pi^j_k(t_1^{j},t_2^{j},t_3^{j},\ldots)=t_k^{j}
$$
for each $(t_1^{j},t_2^{j},t_3^{j},\ldots)\in Q_j$ and 
$$
{\color{blue}p^j_k}: P\rightarrow I_k^j
$$
to denote the composition
$$
p_k^j=\pi_k^{j}\circ p_j.
$$
First, we show that $F$ is compact.  Let $(\mathbf x_k)$ be a sequence of points 
in $F$ and let $\mathbf x\in P$ be such a point that 
$$
\lim_{k\to\infty}\mathbf x_k=\mathbf x.
$$
We show that $F$ is compact by showing that $\mathbf x\in F$.
Note that for all positive integers $j$, $k$ and $n$,
\begin{enumerate}
\item $p_j(\mathbf x_k)=\sigma_j(p_{j+1}(\mathbf x_k))$, 
\item $\displaystyle \lim_{k\to\infty}p_n^j(\mathbf x_k)=p_n^j(\mathbf x)$ and
\item $p_{n+1}^j(\mathbf x_k)$ is equal either to $r\cdot p_n^j(\mathbf x_k)$ or $\rho\cdot p_n^j(\mathbf x_k)$. 
\end{enumerate}
We show that $\mathbf x\in F$ by showing that for all positive integers $j$ and $n$,
\begin{enumerate}
\item[(a)] $p_j(\mathbf x)=\sigma_j(p_{j+1}(\mathbf x))$,  and
\item[(b)] $p_{n+1}^j(\mathbf x)$ is equal either to $r\cdot p_n^j(\mathbf x)$ or $\rho\cdot p_n^j(\mathbf x)$. 
\end{enumerate}
 {To show (a),} let $j$ be any positive integer.  Then 
$$
p_j(\mathbf x)=p_j(\lim_{k\to\infty}\mathbf x_k)=\lim_{k\to\infty}p_j(\mathbf x_k)=\lim_{k\to\infty}\sigma_j(p_{j+1}(\mathbf x_k))=\sigma_j(p_{j+1}(\lim_{k\to\infty}\mathbf x_k))=\sigma_j(p_{j+1}(\mathbf x)).
$$
{To show (b),} let $j$ and $n$ be positive integer.  Then we consider the following possible cases.
\begin{enumerate}
\item There is a positive integer $k_0$ such that for each $k\geq k_0$,  $p_{n+1}^j(\mathbf x_k)=\rho\cdot p_n^j(\mathbf x_k)$.   In this case, 
$$
p_{n+1}^j(\mathbf x)=p_{n+1}^j(\lim_{k\to\infty}\mathbf x_k)=\lim_{k\to\infty}p_{n+1}^j(\mathbf x_k)=\lim_{k\to\infty}\rho\cdot p_n^j(\mathbf x_k)=\rho\cdot p_n^j(\lim_{k\to\infty}\mathbf x_k)=\rho\cdot p_n^j(\mathbf x).
$$
\item There is a positive integer $k_0$ such that for each $k\geq k_0$,  $p_{n+1}^j(\mathbf x_k)=r\cdot p_n^j(\mathbf x_k)$.    In this case, 
$$
p_{n+1}^j(\mathbf x)=p_{n+1}^j(\lim_{k\to\infty}\mathbf x_k)=\lim_{k\to\infty}p_{n+1}^j(\mathbf x_k)=\lim_{k\to\infty}r\cdot p_n^j(\mathbf x_k)=r\cdot p_n^j(\lim_{k\to\infty}\mathbf x_k)=r\cdot p_n^j(\mathbf x).
$$
\item There are strictly increasing sequences $(i_k)$ and $(j_k)$ of positive integers such that for each positive integer $k$, $p_{n+1}^j(\mathbf x_{i_k})=\rho\cdot p_n^j(\mathbf x_{i_k})$ and $p_{n+1}^j(\mathbf x_{j_k})=r\cdot p_n^j(\mathbf x_{j_k})$. In this case, 
\begin{align*}
p_{n+1}^j(\mathbf x)=&p_{n+1}^j(\lim_{k\to\infty}\mathbf x_k)=p_{n+1}^j(\lim_{k\to\infty}\mathbf x_{i_k})=\lim_{k\to\infty}p_{n+1}^j(\mathbf x_{i_k})=\lim_{k\to\infty}\rho\cdot p_n^j(\mathbf x_{i_k})=
\\
&\rho\cdot p_n^j(\lim_{k\to\infty}\mathbf x_{i_k})=\rho\cdot p_n^j(\lim_{k\to\infty}\mathbf x_{k})=\rho\cdot p_n^j(\mathbf x)
\end{align*}
and
\begin{align*}
p_{n+1}^j(\mathbf x)=&p_{n+1}^j(\lim_{k\to\infty}\mathbf x_k)=p_{n+1}^j(\lim_{k\to\infty}\mathbf x_{j_k})=\lim_{k\to\infty}p_{n+1}^j(\mathbf x_{j_k})=\lim_{k\to\infty}r\cdot p_n^j(\mathbf x_{j_k})=
\\
&r\cdot p_n^j(\lim_{k\to\infty}\mathbf x_{j_k})=r\cdot p_n^j(\lim_{k\to\infty}\mathbf x_{k})=r\cdot p_n^j(\mathbf x).
\end{align*}
It follows from $\rho\cdot p_n^j(\mathbf x) = r \cdot p_n^j(\mathbf x)$ and $\rho\neq r$ that $p_n^j(\mathbf x)=0$. Therefore, $p_{n+1}^j(\mathbf x)=0$,  and $p_{n+1}^j(\mathbf x)=\rho\cdot p_n^j(\mathbf x)$ and $p_{n+1}^j(\mathbf x)=r\cdot p_n^j(\mathbf x)$ trivially follow.
\end{enumerate}
It follows that $\mathbf x\in F$. Therefore, $F$ is compact. 

Next, note that for each $(\mathbf a,\mathbf c)\in \{r,\rho\}^{\mathbb N}\times \{r,\rho\}^{\mathbb N}$, $A_{\mathbf a,\mathbf c}$ is the convex arc in $\prod_{j=1}^{\infty}Q_j$ with end-points {$\mathbf O$}  and $\Big({K_1^{\mathbf a,\mathbf c}(1),K_2^{\mathbf a,\mathbf c}(1),K_3^{\mathbf a,\mathbf c}(1)},\ldots\Big)$. 
 For each $(\mathbf a,\mathbf c)\in \{r,\rho\}^{\mathbb N}\times \{r,\rho\}^{\mathbb N}$,  let
 $$
 \mathbf e_{\mathbf a,\mathbf c}=\Big({K_1^{\mathbf a,\mathbf c}(1),K_2^{\mathbf a,\mathbf c}(1),K_3^{\mathbf a,\mathbf c}(1)},\ldots\Big)
$$
and let
$$
E=\{\mathbf e_{\mathbf a,\mathbf c}  \ | \ (\mathbf a,\mathbf c)\in \{r,\rho\}^{\mathbb N}\times \{r,\rho\}^{\mathbb N}\}.
$$
Let $\varphi:\{r,\rho\}^{\mathbb N}\times \{r,\rho\}^{\mathbb N}\rightarrow E$  be defined by
$$
\varphi({\mathbf a,\mathbf c})=\Big({K_1^{\mathbf a,\mathbf c}(1),K_2^{\mathbf a,\mathbf c}(1),K_3^{\mathbf a,\mathbf c}(1)},\ldots\Big)
$$
for each $({\mathbf a,\mathbf c})\in \{r,\rho\}^{\mathbb N}\times  \{r,\rho\}^{\mathbb N}$. Obviously, $\varphi$ is a homeomorphism from the Cantor set $\{r,\rho\}^{\mathbb N}\times  \{r,\rho\}^{\mathbb N}$ to $E$.   Therefore, $E$ is a Cantor set and it follows that $F$ is a Cantor fan with $E(F)=E$.  Since $M$ is a subcontinuum of $F$, it is itself a smooth fan.  
\end{proof}
\begin{lemma}\label{tazadnjalema}
Let $\mathbf x\in M$.  If for each positive integer $j$, $p_j(\mathbf x)\in E(M_{r,\rho})$, then  $\mathbf x\in E(M)$.
\end{lemma}
\begin{proof}
Let $t\in(0,1]$ and let $(\mathbf a,\mathbf c)\in \mathcal P_{r,\rho}$ be such that 
$$
\mathbf x=\Big({K_1^{\mathbf a,\mathbf c}(t),K_2^{\mathbf a,\mathbf c}(t),K_3^{\mathbf a,\mathbf c}(t)},\ldots\Big). 
$$
Suppose that for each positive integer $j$, $p_j(\mathbf x)\in E(M_{r,\rho})$.  It follows from Theorem \ref{krajisca} that for each positive integer $n$,
\begin{align*}
\sup\Big\{&\frac{t}{c_{n-1}c_{n-2}c_{n-3}\ldots c_3c_2c_1}, \frac{t}{c_{n-2}c_{n-3}c_{n-4}\ldots c_3c_2c_1}, \frac{t}{c_{n-3}c_{n-4}c_{n-5}\ldots c_3c_2c_1}, \ldots ,
\\
&\frac{t}{c_3c_2c_1}, \frac{t}{c_2c_1}, \frac{t}{c_1}, t,a_1\cdot t, a_2a_1\cdot t, a_3a_2a_1\cdot t,\ldots\Big\}= 1.
\end{align*}
Then $\mathbf x$ is an end point of the arc $L_{(\mathbf a,\mathbf c)}$ and it follows that $\mathbf x\in E(M)$.
\end{proof}
\begin{theorem}\label{mainresult}
$M$ is a Lelek fan.
\end{theorem}
\begin{proof}
It follows from Theorem \ref{fencic} that $M$ is a smooth fan.  To complete the proof, we need to show that $E(M)$ is dense in $M$.  To see that, let $\mathbf t\in M$
 be any point and let $\varepsilon>0$ such that $\varepsilon <1$. We show that $B(\mathbf t,\varepsilon)\cap E(M)\neq \emptyset$.  Since $\mathbf t\in M$, there is a useful pair $(\mathbf a,\mathbf c)\in  \mathcal P_{r,\rho}$ such that 
 $$
 \mathbf t=\Big({K_1^{\mathbf a,\mathbf c}(t),K_2^{\mathbf a,\mathbf c}(t),K_3^{\mathbf a,\mathbf c}(t)},\ldots\Big)
 $$
 for some $t\in [0,1]$. Fix such an $(\mathbf a,\mathbf c)=\Big((a_1,a_2,a_3,\ldots),(c_1,c_2,c_3,\ldots)\Big)\in  \mathcal P_{r,\rho}$ and $t\in [0,1]$. 
 If $t=1$, then $\mathbf t\in E(M)$ and there is nothing to prove.  For the rest of the proof, we assume that $t<1$.  We consider the following possible cases.
 \begin{enumerate}
   \item $t=0$.  It follows that $\mathbf t=\mathbf O$.  For each positive integer $n$, let $c_n=\rho$, and let $\mathbf e=(e_1,e_2,e_3,\ldots)\in E(M_{r,\rho})$ { be any point}.  For each positive integer $n$, let
    \begin{align*}
  \mathbf e_n=\Bigg(&\Big(\frac{e_1}{c_nc_{n-1}c_{n-2}\ldots c_3c_2c_1},\ldots,\frac{e_1}{c_3c_2c_1},\frac{e_1}{c_2c_1},\frac{e_1}{c_1},e_1,e_1,e_2,e_3,\ldots\Big),\\
&\Big(\frac{e_1}{c_{n+1}c_{n}c_{n-1}\ldots c_3c_2c_1},\ldots,\frac{e_1}{c_3c_2c_1},\frac{e_1}{c_2c_1},\frac{e_1}{c_1},e_1,e_1,e_2,e_3,\ldots\Big),\\
&\Big(\frac{e_1}{c_{n+2}c_{n+1}c_{n}\ldots c_3c_2c_1},\ldots,\frac{e_1}{c_3c_2c_1},\frac{e_1}{c_2c_1},\frac{e_1}{c_1},e_1,e_1,e_2,e_3,\ldots\Big),\ldots\Bigg).
 \end{align*}

%  \begin{align*}
%	W_n &= O'_n, \\
%	W_{n-1} &=  O'_{n-1} \cap  \frac{1}{a_{n-1}} W_{n}=  O'_{n-1} \cap  \frac{1}{a_{n-1}} O'_{n}, \\
%	W_{n-2} &=  O'_{n-2} \cap  \frac{1}{a_{n-2}} W_{n-1}=  O'_{n-2} \cap  \frac{1}{a_{n-2}} O'_{n-1}\cap \frac{1}{a_{n-2}a_{n-1}} O'_{n}, \\
%	W_{n-3} &=  O'_{n-3} \cap  \frac{1}{a_{n-3}} W_{n-2}=  O'_{n-3} \cap  \frac{1}{a_{n-3}} O'_{n-2}\cap \frac{1}{a_{n-3}a_{n-2}} O'_{n-1}\cap   \frac{1}{a_{n-3}a_{n-2}a_{n-1}} O'_{n},  \\
%	&\vdots\\
%	W_2 &=O'_{2} \cap  \frac{1}{a_{2}} W_{3}=O'_{2} \cap \frac{1}{a_2} O'_{3} \cap \frac{1}{a_2 a_3} O'_{4} \cap \ldots \cap \frac{1}{a_2a_3a_4 \ldots a_{n-1}} O'_{n}, \\
%	W_1 &=O'_{1} \cap  \frac{1}{a_{1}} W_{2}=O'_{1} \cap \frac{1}{a_1} O'_{2} \cap \frac{1}{a_1 a_2} O'_{3} \cap  \frac{1}{a_1 a_2a_3} O'_{4} \cap \ldots \cap \frac{1}{a_1a_2a_3 \ldots a_{n-1}} O'_{n}.
%	\end{align*}
%	
%	
	
  Obviously, for each positive integer $n$, $ \mathbf e_n\in E(M)$ and $\displaystyle \lim_{n\to \infty}\mathbf e_n=\mathbf O=\mathbf t$.
 \item $t\neq 0$.  Note that since  $M_{r,\rho}$ is a Lelek fan, the set $E(M_{r,\rho})$ is dense in $M_{r,\rho}$. Therefore,  
 $$
 B\Big((t,a_1\cdot t, a_2a_1\cdot t, a_3a_2a_1\cdot t,\ldots),\frac{\varepsilon}{2}\Big)\cap E(M_{r,\rho})\neq \emptyset.
 $$  
 Also, note that since $(\mathbf a,\mathbf c)$ is a useful pair, it follows that the sequence
 $$
\left( \frac{1}{c_nc_{n-1}c_{n-2}\ldots c_3c_2c_1 }\right)
 $$
 is bounded from above.  Thus,  $\inf\{c_nc_{n-1}c_{n-2}\ldots c_3c_2c_1 \ | \ n \textup{ is a positive integer}\}$ does exist and it is not equal to $0$, i.e., it is a positive number. 
 Therefore, there is an end point 
 $$
 \mathbf e=(e_1,e_2,e_3,\ldots)\in E(M_{r,\rho})
 $$
 of $M_{r,\rho}$ such that {(a), (b) and (c) below hold:}
 \begin{itemize}
 \item[(a)] $ e_1\in (0,t)$,
 \item[(b)] $ t-e_1<\inf\{c_nc_{n-1}c_{n-2}\ldots c_3c_2c_1 \ | \ n \textup{ is a positive integer}\}\cdot \frac{\varepsilon}{2}$,  and 
 \item[(c)] $d\Big(\mathbf e,(t,a_1\cdot t, a_2a_1\cdot t, a_3a_2a_1\cdot t,\ldots)\Big)<\frac{\varepsilon}{2}$. 
 \end{itemize}
{Chose and fix such an endpoint $\mathbf e$.}  It follows from Theorem \ref{krajisca} that $\sup\{e_1,e_2,e_3,\ldots\}=1$. 
 Let 
\begin{align*}
\mathbf x=\Bigg(&\Big(e_1,e_2,e_3,\ldots\Big),
\Big(\frac{e_1}{c_1},e_1,e_2,e_3,\ldots \Big),
\Big(\frac{e_1}{c_2c_1},\frac{e_1}{c_1},e_1,e_1,e_2,e_3,\ldots\Big),\\
&\Big(\frac{e_1}{c_3c_2c_1},\frac{e_1}{c_2c_1},\frac{e_1}{c_1},e_1,e_1,e_2,e_3,\ldots\Big),\ldots\Bigg).
\end{align*}
  Note that 
\begin{align*}
&\sup\Big\{\frac{e_1}{c_nc_{n-1}c_{n-2}\ldots c_3c_2c_1} \ | \ n \textup{ is a positive integer}\Big\}\leq 
\\
& \sup\Big\{\frac{t}{c_nc_{n-1}c_{n-2}\ldots c_3c_2c_1} \ | \  n \textup{ is a positive integer}\Big\}\leq 1.
\end{align*}
 It follows that $\mathbf x\in M$. Since $\sup \{e_1,e_2,e_3,\ldots\}=1$, it follows that for each positive integer $n$,
 $$
 \sup\Big\{\frac{e_1}{c_nc_{n-1}c_{n-2}\ldots c_3c_2c_1},\ldots,\frac{e_1}{c_3c_2c_1},\frac{e_1}{c_2c_1},\frac{e_1}{c_1},e_1,e_1,e_2,e_3,\ldots\Big\}=1.
 $$
 It follows from Theorem \ref{krajisca} that for each positive integer $n$, $p_n(\mathbf x)\in E(M_{r,\rho})$. 
 By Lemma \ref{tazadnjalema}, $\mathbf x\in E(M)$.  
 
 Next, we prove that 
 $$
 D(\mathbf t,\mathbf x)=\max\left\{\frac{d(p_n(\mathbf x),p_n(\mathbf t))}{2^n} \ | \ n \textup{ is a positive integer}\right\}<\varepsilon. 
 $$
 We show this by showing that for each positive integer $n$, $ d(p_n(\mathbf x),p_n(\mathbf t))<\frac{\varepsilon}{2}$.
Let $n$ be a positive integer.  Note that 
$$
d\Big(\mathbf e,(t,a_1\cdot t, a_2a_1\cdot t, a_3a_2a_1\cdot t,\ldots)\Big)<\frac{\varepsilon}{2}.
$$
 Also, note  that for each $j\in \{1,2,3,\ldots,n\}$,
\begin{align*}
 &\frac{t}{c_jc_{j-1}c_{j-2}\ldots c_3c_2c_1}-\frac{e_1}{c_jc_{j-1}c_{j-2}\ldots c_3c_2c_1}=\\
 &\frac{1}{c_jc_{j-1}c_{j-2}\ldots c_3c_2c_1}(t-e_1)<
\\
&\frac{\inf\{c_nc_{n-1}c_{n-2}\ldots c_3c_2c_1 \ | \ n \textup{ is a positive integer}\}\cdot \frac{\varepsilon}{2}}{c_jc_{j-1}c_{j-2}\ldots c_3c_2c_1}\leq\frac{\varepsilon}{2},
\end{align*}
 since $\inf\{c_nc_{n-1}c_{n-2}\ldots c_3c_2c_1 \ | \ n \textup{ is a positive integer}\}\leq c_jc_{j-1}c_{j-2}\ldots c_3c_2c_1$.
 It follows that $ d(p_n(\mathbf x),p_n(\mathbf t))<\frac{\varepsilon}{2}$ and, therefore,  $D(\mathbf t,\mathbf x)<\varepsilon$.
 \end{enumerate}
 This completes the proof. 
\end{proof}

We conclude the paper by stating two open problems. 
 \begin{problem}\label{pr0}
Find another path-connected non locally connected continuum (not homeomorphic to the Lelek fan)  that  admits a transitive homeomorphism.
\end{problem}

We need the following definition in the statement of Problem \ref{pr1}. The definition defines a property, which is stronger than the topological transitivity.
\begin{definition}
Let $X$ be a continuum and let $f:X\rightarrow X$ be a continuous function. We say that $f$ is \emph{\color{blue} topologically mixing}, if for any two non-empty open sets $U$ and $V$ in $X$, there is a positive integer $n_0$ such that for all integers $n\geq n_0$, $f^n(U)\cap V\neq \emptyset$.
\end{definition}
%\begin{problem}\label{pr1}
%Is the shift map $\sigma$ on the Lelek fan $M$  topologically mixing?
%\end{problem}
\begin{problem}\label{pr1}
Is there a topologically mixing homeomorphism on the Lelek fan? 
\end{problem}

\section{Acknowledgement}
This work is supported in part by the Slovenian Research Agency (research project J1-4632 and research program P1-0285).
	
%\section{Declarations}
%The following sections are not relevant to our manuscript. 
%\subsection{Competing interests}
%Not applicable.
%\subsection{Data Availability Statement}
%Not applicable.

\noindent I. Bani\v c\\
              (1) Faculty of Natural Sciences and Mathematics, University of Maribor, Koro\v{s}ka 160, SI-2000 Maribor,
   Slovenia; \\(2) Institute of Mathematics, Physics and Mechanics, Jadranska 19, SI-1000 Ljubljana, 
   Slovenia; \\(3) Andrej Maru\v si\v c Institute, University of Primorska, Muzejski trg 2, SI-6000 Koper,
   Slovenia\\
             {iztok.banic@um.si}           %  \\
     
				\-
				
		\noindent G.  Erceg\\
             Faculty of Science, University of Split, Rudera Bo\v skovi\' ca 33, Split,  Croatia\\
%             {i}     
{{gorerc@pmfst.hr}       }    %  \\

                 \-

                 	\-
					
  \noindent J.  Kennedy\\
             Lamar University, 200 Lucas Building, P.O. Box 10047, Beaumont, TX 77710 USA\\
%             {}     
{{kennedy9905@gmail.com}       }    %  \\

%``text''
%%%%%%%%%%%%%%%%%%%%%%%%%%%%%%%%%%%%%%%%%%%%%%%%%%%%%%%%%%%%%%%%%%%%%%%%%%%%%%%%%
%%% I N T R O D U C T I O N S


\begin{thebibliography}{9}
\bibitem{A} E. ~Akin, {General Topology of Dynamical Systems}, Volume 1, Graduate Studies in Mathematics Series, American Mathematical Society, Providence RI, 1993.
\bibitem{banic1} I.~Bani\v c, G.~Erceg,  J.~Kennedy, The Lelek fan as the inverse limit of intervals with a single set-valued bonding function whose graph is an arc, 
https://doi.org/10.48550/arXiv.2206.00087
\bibitem{barge} M.~Barge, J.~Martin, Chaos, periodicity, and snakelike continua. Trans. Amer. Math. Soc. 289 (1985), no. 1, 355--365.
\bibitem{BK} D.~Barto\v sova,  A.~Kwiatkowska, Lelek fan from a projective Fraisse limit, Fundamenta Math.  231 (2015) 57--79.
\bibitem{jan} J. ~Boro\' nski, P. ~Minc and S.~ \v Stimac,  On conjugacy between natural extensions of 1-dimensional maps,  Ergod. Th.  Dynam. Sys.   (2022) https://doi.org/10.1017/etds.2022.62.
\bibitem{jan2} J.~ Boro\' nski,  J. ~Kupka, New chaotic planar attractors from smooth zero entropy interval maps, Adv. Difference Equ.  232 (2015) 11 pp.
\bibitem{jan3} J.~ Boro\' nski,  P.~Oprocha, On indecomposability in chaotic attractors. Proc. Amer. Math. Soc. 143 (2015),  3659--3670.
\bibitem{jan4} J.~ Boro\' nski,  P.~Oprocha, On dynamics of the Sierpin\' ski carpet,  C. R. Math. Acad. Sci. Paris 356 (2018) 340--344.
\bibitem{oversteegen} W.~D.~Bula and L.~Overseegen, A Characterization of smooth Cantor Bouquets,  Proc. Amer.Math.Soc. 108 (1990) 529--534.
\bibitem{charatonik} W.~J.~Charatonik, The Lelek fan is unique, Houston J. Math. 15 (1989) 27--34.
\bibitem{charatonik2} J. ~J.~ Charatonik, W. ~J. ~Charatonik and S. ~Miklos, Confluent mappings of fans. Dissertationes Math. (Rozprawy Mat.) 301 (1990), 86 pp.
\bibitem{cinc} J.~\v Cin\v c, P.~Oprocha, Parametrized family of pseudo-arc attractors: Physical measures and prime end rotations, Proc. London Math. Soc.  125 (2022) 318--357.
\bibitem{handel} M. ~Handel, A pathological area preserving $C^{\infty}$ diffeomorphism of the plane, Proc. Amer.Math.Soc.86 (1982),163--168
\bibitem{he} F. ~He, J. ~Liu, Invariant measures and uniform positive entropy property for inverse limits, Appl. Math. J. Chinese Univ. Ser. B. 14 (1999)  265--272.
\bibitem{hernandez} R.~Hernandez - Gutierrez, L.~C.~Hoehn, Smooth fans that are endpoint rigid, 
https://doi.org/10.48550/arXiv.2206.12776.
\bibitem{HM} L.~C. ~Hoehn and C. ~Mouron, Hierarchies of chaotic maps on continua, Ergodic Theory Dynam. Systems 34 (2014), 1897--1913.
\bibitem{judy} J. ~Kennedy,   A transitive homeomorphism on the pseudoarc which is semiconjugate to the tent map, Trans. Amer. Math. Soc. 326 (1991),  773--793.
\bibitem{KS} S.~Kolyada, L.~Snoha,  Topological transitivity, \textit{Scholarpedia}  4 (2):5802 (2009).
\bibitem{lelek} A.~Lelek, On plane dendroids and their end-points in the classical sense, Fund. Math. 49 (1960/1961) 301--319.
\bibitem{li} S. ~Li, Dynamical properties of the shift maps on the inverse limit spaces,  Ergod. Th.  Dynam. Sys.  12 (1992) 95--108.
\bibitem{chris1} V. ~Mart\' nez-de-la-Vega, J.~M. ~Mart\' inez-Montejano, C.~Mouron, Mixing homeomorphisms and indecomposability,  Topology App. \textbf{254} (2019) 50--58.
\bibitem{minc} P. ~Minc and W. ~R.~ R. ~Transue, A Transitive Map on [0,1] Whose Inverse Limit is the Pseudoarc, Proceedings of the American Mathematical Society 111 (1991) 1165--1170. 
\bibitem{chris2} C.~Mouron, Tree-like continua do not admit expansive homeomorphisms. Proceedings of the A.M.S. 130 Nov. 2002, p. 3409-3413.
\bibitem{chris3} C.~Mouron,  Positive entropy homeomorphisms of chainable continua and indecomposable subcontinua, Proc. Amer. Math. Soc. 139 (2011), no. 8, 2783--2791.
\bibitem{chris4} C.~Mouron,  Expansive homeomorphisms and indecomposable subcontinua. Topology Appl. 126 (2002), no. 1-2, 13--28.
\bibitem{chris5} C.~Mouron,  Mixing sets, positive entropy homeomorphisms and non-Suslinean
continua. Ergodic Theory and Dynamical Systems, 36 (2016), no. 7, 2246--2257.
\bibitem{nadler} S.~B.~Nadler, Continuum theory. An introduction, Marcel Dekker, Inc., New York, 1992.
%\bibitem{perez} R. ~Perez-Marco, Hedgehog’s Dynamics, Preprint, UCLA, 1996.
\bibitem{seidler} G.~T.~ Seidler, The topological entropy of homeomorphisms on one-dimensional continua,Proc. Amer. Math. Soc. 4 (1990), 1025--1030.
\end{thebibliography}
\end{document}